\documentclass[12pt]{article}
\usepackage{amsfonts}
\usepackage[dvips]{graphics}

\usepackage[utf8]{inputenc}
\usepackage[english]{babel}
\usepackage{latexsym,cite,mathrsfs}
\usepackage{amsfonts}
\usepackage{amssymb, amsthm}
\usepackage{xcolor}
\usepackage{marginnote}
\usepackage{amsmath}
\usepackage{a4wide}
\usepackage[affil-it]{authblk}

\newtheorem{theorem}{Theorem}

\newtheorem{definition}{Definition}
\newtheorem{remark}{Remark}
\newtheorem{example}{Example}

\newtheorem{proposition}{Proposition}
\begin{document}
\title {A Liouville comparison principle for solutions of
semilinear elliptic  second-order partial differential inequalities}
\author{Vasilii V. Kurta}
\maketitle
\begin{abstract}
\noindent We consider semilinear elliptic second-order partial
differential inequalities of the form
$$
Lu + |u|^{q-1}u \leq  Lv + |v|^{q-1}v \qquad (\ast)
$$
in the whole space ${\mathbb R}^n$,  where $n\geq 2$, $q>0$ and $L$
is a linear elliptic  second-order partial differential operator in
divergence form. We assume that the coefficients  of the operator
$L$ are measurable and  locally bounded such that the quadratic form associated with the operator $L$ is
symmetric and  non-negative definite. We obtain a Liouville
comparison principle in terms of a capacities associated with the
operator $L$ for solutions of ($\ast$) which are 
measurable   and belong locally in  ${\mathbb R}^n$ to a
Sobolev-type function space also  associated  with the operator $L$.
\end{abstract}

\thispagestyle{empty}
\section{Definitions}
Let $L$ be a linear elliptic second-order  partial differential
operator defined by 
\begin{eqnarray}
\label{(1)}
Lu:=\sum\limits_{i,j=1}^n (a_{ij}(x)u_{x_i})_{x_j}
\end{eqnarray}
in the whole space  ${\mathbb R}^n$, $n\geq2$. The
coefficients $a_{ij}$ of the  operator $L$ are measurable and locally bounded in
${\mathbb R}^n$, $a_{ij}=a_{ji}$ , $i,j=1, \dots, n$. We assume that the
quadratic form associated with the operator $L$ is non-negative definite, namely, it satisfies the condition
\begin{eqnarray}\label{(3)}
\sum\limits_{i,j=1}^n a_{ij}(x)\xi_{i}\xi_{j}\geq 0
\end{eqnarray}
for all $\xi=(\xi_1, \dots, \xi_n)\in {\mathbb R}^n$ and almost all
$x\in {\mathbb R}^n$.

We study  solutions $(u,v)$ of semilinear elliptic second-order
partial differential inequalities of the form
\begin{equation}\label{(4)}
{L}u + |u|^{q-1}u \leq {L}v + |v|^{q-1}v,
\end{equation}
where $q>0$.
We assume that  $(u,v)$ is measurable in ${\mathbb
R}^n$ and belong locally to a Sobolev-type function space
associated with the differential
operator $L$.
We  call such solution the entire solution  of
the inequality \eqref{(4)} in ${\mathbb R}^n$.

Note that if $u$ and $v$ satisfy, respectively, the semilinear
elliptic second-order partial differential inequalities
\begin{equation}\label{(5)}
- {L} u\geq|u|^{q-1}u 
\end{equation}
and
\begin{equation}\label{(6)}
- {L} v\leq |v|^{q-1}v, 
\end{equation}
then the pair $(u, v)$  satisfies the inequality \eqref{(4)}. Hence, all the
results obtained in this paper for solutions of \eqref{(4)} are valid  for
the corresponding solutions of the system \eqref{(5)}--\eqref{(6)}.

\begin{definition} Let $n\geq 2$, $q>0$, $\hat q = \max \{1, q\}$,  $L$ be a
differential operator defined by  \eqref{(1)} in $ {\mathbb
R}^n$, and  $\mathcal Q$ be a bounded domain in ${\mathbb R}^n$.
By $W^{L,q}(\mathcal Q)$ we denote the completion of the function
space $C^{\infty} (\mathcal Q)$ with respect to the norm defined by
the expression
\begin{eqnarray}\label{(7)}
\|f\|_{W^{L,q}(\mathcal Q)}:= \left(\int\limits_{\mathcal
Q}\sum\limits_{i,j=1}^n
a_{ij}f_{x_i}f_{x_j}dx\right)^{1/2}+\left(\int\limits_{\mathcal
Q}|f|^{\hat q}dx\right)^{1/{\hat q}}.
\end{eqnarray}
\end{definition}

\begin{definition} Let $n\geq 2$, $q>0$,  and  $L$ be a  differential operator defined by  \eqref{(1)} in $ {\mathbb R}^n$. We say that a  function $f$
belongs to the function space $W^{L,q}_{\mathrm {loc}}({\mathbb
R}^n)$,  if it belongs to  $W^{L,q}(\mathcal Q)$ for every bounded
domain $\mathcal Q$ in ${\mathbb R}^n$.
\end{definition}

\begin{definition} Let $n\geq 2$, $q>0$,  and   $L$ be a  differential operator defined by  \eqref{(1)} in $ {\mathbb R}^n$.
We say that a pair of functions $(u,v)$ is  an  entire solution of the inequality \eqref{(4)}  in
${\mathbb R}^n$,  if $u$ and $v$  belong to the function space
$W^{L,q}_{\mathrm
{loc}}({\mathbb R}^n)$ and satisfy the integral inequality
\begin{eqnarray}\label{(8)}
\int\limits_{{\mathbb R}^n}\left( \sum\limits_{i,j=1}^n
a_{ij} u_{x_i}\zeta_{x_j}-|u|^{q-1}u\zeta\right) dx \geq
\int\limits_{{\mathbb R}^n}\left(\sum\limits_{i,j=1}^n
a_{ij}v_{x_i}\zeta_{x_j}-|v|^{q-1}v\zeta\right)dx
\end{eqnarray}
for every non-negative function  $\zeta\in C^{\infty}({\mathbb
R}^n)$ with compact support.
\end{definition}

We understand  inequality \eqref{(8)}  in the sense that discussed as in
\cite{[4]} or \cite{[9]}.
Definitions of entire solutions to  inequalities \eqref{(5)}--\eqref{(6)},
are the special cases of Definition 3 with
$v\equiv 0$ or $u\equiv 0$, respectively.

\begin{definition}
Let $n\geq2$.
We denote by ${\rm{Lip}}_{\rm{loc}}({\mathbb R}^n)$  the space of measurable functions $f$ defined in ${\mathbb R}^n$  which satisfy the Lipschitz
condition
\begin{eqnarray*}\label{(9)}
|f(x)-f(y)|\leq K|x-y|
\end{eqnarray*}
on any compact set $\mathcal K$  in ${\mathbb R}^n$, with $K$ some
positive constant which possibly depends on $\mathcal K$.
\end{definition}

\begin{definition}
Let $n\geq2$,  $\mathcal E$ be a domain in ${\mathbb R}^n$, and
$\mathcal D$, $\mathcal G$ be two sets in $\mathcal E$, which
are non-overlapping and closed with respect to $\mathcal E$. We call
such a triple $(\mathcal D,\mathcal G;\mathcal E)$ the condenser.
\end{definition}

For a given differential  operator $L$   defined by  \eqref{(1)} in
$ {\mathbb R}^n$, following \cite{[7]}, we use  a notion of the
$(L,p)$-capacity of the condenser $(\mathcal D,\mathcal G;\mathcal
E)$ in the particular case when $\mathcal E$ coincides with the
whole space ${\mathbb R}^n$.

\begin{definition}  Let $n\geq 2$, $p> 1$,  and   $L$ be a  differential operator defined by the
relation \eqref{(1)} in $ {\mathbb R}^n$. We call
the quantity
\begin{eqnarray}\label{(10)}
{\mathrm {cap}}_{L,p} (\mathcal D,\mathcal G;{\mathbb R}^n)= \inf
\left [\int\limits_{{\mathbb R}^n}\left(\sum\limits_{i,j=1}^n
a_{ij}{\varphi}_{x_i}{\varphi}_{x_j}\right)^{p/2}dx\right]
\end{eqnarray}
the $(L,p)$-capacity of the condenser  $(\mathcal D,\mathcal
G;{\mathbb R}^n)$. Here the infimum is taken over the all functions
$\varphi$ in the space ${\rm{Lip}}_{\rm{loc}}({\mathbb
R}^n)$ such that $\varphi=1$ on $\mathcal D$, $\varphi=0$ on
$\mathcal G$, and $1\geq \varphi \geq 0$ in  $ {\mathbb R}^n$.
\end{definition}

In the case when the coefficients $a_{ij}$ of a  differential
operator $L$ coincide with  Kronecker's  symbols
\begin{eqnarray*}
\delta_{ij}=\left \{ \begin {array} {ll} 1, &
{\rm if}\   i=j,\\
0, & {\rm if}\  i\neq j, \end{array}\right.
\end{eqnarray*}
the $(L,p)$-capacity ${\mathrm {cap}}_{L,p}$ associated
with the operator $L$ is denoted by ${\mathrm {cap}}_{p}$. Notice that if $p=2$, then ${\mathrm {cap}}_{2}$ is the Wiener capacity. 

Now, we give an estimate for the  $(L,p)$-capacity of the
condenser $(\overline {\mathcal B}_r,{\mathbb R}^n \setminus
{\mathcal B}_R ;{\mathbb R}^n)$ via the well-known nonlinear
$p$-capacity ${\mathrm {cap}}_p$ of the same condenser,  and with the
coefficients $a_{ij}(x)$ of the operator $L$ defined on the set ${{\mathcal
B}_R\setminus {\mathcal B}_r}$ for any $R>r>1$, where ${\mathcal
B}_r=\{x\in {\mathbb R}^n: |x|< r\}$ and $\overline {\mathcal
B}_r={\{x\in {\mathbb R}^n: |x|\leq r\}}$.

\begin{proposition}\label{prop1} Let $n\geq 2$, $p > 1$,  and let  $L$ be a  differential operator defined by  \eqref{(1)} in $ {\mathbb R}^n$. Then
\begin{eqnarray}\label{(12)}
&&{\mathrm {cap}}_{L,p} (\overline {\mathcal B}_r,{\mathbb
R}^n\setminus {\mathcal B}_R ;{\mathbb R}^n)\leq\nonumber  \\
&&~~~~~~~~~~~~~~~~\sup_{x\in {\mathcal B}_R\setminus {\mathcal B}_r}
\left(\sum\limits_{i,j=1}^n a^2_{ij}(x)\right)^{p/4} {\mathrm
{cap}}_p (\overline {\mathcal B}_r,{\mathbb R}^n \setminus {\mathcal
B}_R ;{\mathbb R}^n)
\end{eqnarray}
for all $R>r>1$.
\end{proposition}

It is well known (see, e.g., \cite{[2]}, p.178 or  \cite{[8]}, p.12)
that for any $n\geq 2$ and $p>1$ the inequality
\begin{eqnarray}\label{(13)}
 {\mathrm {cap}}_p ({\overline {\mathcal B}}_{R/2},{\mathbb R}^n\setminus {\mathcal B}_R
;{\mathbb R}^n)\leq C(1+R^2)^{(n-p)/2}
\end{eqnarray}
holds for all $R>0$  with $C$ being some  positive constant which
depends only on $n$ and $p$. Now, based  on  inequalities \eqref{(12)}
and \eqref{(13)}, we give an estimate for the $(L,p)$-capacity of the
condenser $({\overline {\mathcal B}}_{R/2},{\mathbb R}^n \setminus
{\mathcal B}_R ;{\mathbb R}^n)$ in the case when the coefficients
$a_{ij}$ of a differential operator 
\eqref{(1)} satisfy 
\begin{eqnarray}\label{(14)}
\sup_{x\in {\mathcal B}_R\setminus {\mathcal B}_{R/2}}
\sum\limits_{i,j=1}^n a^2_{ij}(x) \leq AR^{-2\sigma}
\end{eqnarray}
for  all sufficiently large $R$, with some
constants $A>0$ and  $\sigma$.

\begin{proposition}\label{prop2} Let $n\geq 2$, $p > 1$,  and let  $L$ be a  differential operator
defined by \eqref{(1)} in $ {\mathbb R}^n$ such that its
coefficients  satisfy condition  \eqref{(14)} for all sufficiently large
$R$, with  some constants $A>0$ and $\sigma$. Then
\begin{eqnarray}\label{(15)}
{\mathrm {cap}}_{L,p} (\overline {\mathcal B}_{R/2},{\mathbb
R}^n\setminus {\mathcal B}_R ;{\mathbb R}^n)\leq {\hat C}R^{(2n-
p(\sigma+2))/2}
\end{eqnarray}
for all sufficiently large $R$, where $\sigma$ is the same constant as
in condition \eqref{(14)} and $\hat C$ is a positive constant which
depends only on $A$, $n$,  $p$ and $\sigma$.
\end{proposition}

\section{Results}

\vskip 5pt \noindent 

For a given differential  operator $L$   defined by
\eqref{(1)} in $ {\mathbb R}^n$ and $q>0$, our goal in this
paper is to establish  a Liouville comparison principle for
solutions of the inequality \eqref{(4)}  defined in  ${{\mathbb R}^n}$ in
terms of the capacity,   associated with the operator $L$, of
condensers of the form  $(\overline{\mathcal B}_{R/2}, {\mathbb
R}^n\setminus {\mathcal B}_R;{\mathbb R}^n)$ as $R\to \infty$. We
illustrate our results in the particular case when the coefficients
of a differential operator $L$ defined by \eqref{(1)} in
${\mathbb R}^n$ satisfy  condition \eqref{(14)}. More precisely, we
obtain  the following results.

\begin{theorem}\label{thm1}
Let $n\geq2$, $q>0$, let  $L$ be a differential operator defined by
\eqref{(1)} in $ {\mathbb R}^n$, let
\begin{eqnarray}\label{(16)}
\liminf\limits_{R\to\infty}{\mathrm {cap}}_{L,2} (\overline
{\mathcal B}_{R/2},{\mathbb R}^n\setminus {\mathcal B}_R; {\mathbb
R}^n )<\infty, \end{eqnarray} 
and let $(u, v)$ be an entire solution
of the inequality \eqref{(4)} in $\mathbb{R}^n$ such that $u(x)\geq v(x)$ a.e. in $\mathbb{R}^n$.
Then $u = v$ a.e. in ${\mathbb R}^n$.
\end{theorem}

Next,  we give a  condition on the coefficients of the operator $L$,
which guarantees that the inequality \eqref{(16)} holds.

\begin{proposition}\label{prop3}
Let $n\geq2$,  and  let  $L$ be a differential operator defined by
 \eqref{(1)} in $ {\mathbb R}^n$  and such that its coefficients
satisfy condition \eqref{(14)} for  all sufficiently large $R$, with some  constants
$A>0$ and  $\sigma \geq n-2$. Then the inequality
\eqref{(16)} holds.
\end{proposition}

The following example shows  the sharpness of Proposition \ref{prop3}.

\begin{example}
Let $n\geq 2$. Consider the differential operator $L$  defined by
the relation
\begin{eqnarray}\label{(17)}
Lu:=\sum\limits_{i=1}^n \frac{\partial }{\partial x_i} \left[\frac
1{(1+|x|^2)^{{\sigma}/ 2}} \frac{\partial u}{\partial x_i} \right]
\end{eqnarray}
in ${\mathbb R}^n$, with some constant $\sigma$. Using standard
arguments in the capacity theory (see, e.g., [2, p.178] or [8, p.12]),
it is not difficult to verify  that the $(L,2)$-capacity, associated
with the  operator L,   of the  condenser $({\overline {\mathcal
B}}_{R/2},{\mathbb R}^n\setminus {\mathcal B}_R; {\mathbb R}^n )$
satisfies the two-sided inequality
\begin{eqnarray}\label{(18)}
C_2 R^{n-\sigma-2}\geq
{\mathrm {cap}}_{L,2} (\overline {\mathcal B}_{R/2},{\mathbb
R}^n\setminus {\mathcal B}_R; {\mathbb R}^n)\geq C_1  R^{n-\sigma-2}
\end{eqnarray}
for all $R>1$, where $\sigma$ is the same constant as in the relation
\eqref{(17)} and $C_1$, $C_2$ some positive constants which depend only on
$n$ and $\sigma$. In turn,  the  inequality \eqref{(18)} yields that the
$(L,2)$-capacity  of the  condenser $(\overline {\mathcal
B}_{R/2},{\mathbb R}^n\setminus {\mathcal B}_R; {\mathbb R}^n )$ is
bounded above by a constant which depends only on $n$ and $\sigma$,
for all $R>1$ when $\sigma\geq n-2$,  and  tends to infinity as
$R\to \infty$ when $\sigma<n-2$.
\end{example}

As a simple corollary of    Theorem \ref{thm1} and Proposition \ref{prop3},  we have  the
following result.

\begin{theorem}\label{thm2}
Let $n\geq2$, $q>0$. Assume that $L$ is a differential operator defined by
\eqref{(1)} in $ {\mathbb R}^n$ such that its coefficients
satisfy condition \eqref{(14)} for  all sufficiently large $R$ with some  constants
$A>0$ and  $\sigma \geq n-2$. Let $(u, v)$ be an
entire solution of the inequality \eqref{(4)} in ${\mathbb R}^n$ such that
$u\geq v$ a.e. in ${\mathbb R}^n$. Then $u = v$ a.e. in ${\mathbb R}^n$.
\end{theorem}

The following example demonstrates the sharpness of Theorem \ref{thm2} in the
case when $1>q>0$.

\begin{example}
Let $n\geq 2$ and  $1>q>0$. Consider the differential operator $L$
from Example 1 defined by \eqref{(17)} with $\sigma<n-2$. It is
easy to verify that the pair $(u,v)$ of the functions
\begin{eqnarray*}\label{(19)}
u(x)=\alpha (1+|x|^2)^{(2+\sigma)/(2(1-q))} + (1+|x|^2)^{-\mu}
\end{eqnarray*}
and
\begin{eqnarray*}\label{(20)}
v(x)=\alpha (1+|x|^2)^{(2+\sigma)/(2(1-q))},
\end{eqnarray*}
where $\alpha$ is a suitable sufficiently large positive constant and
$0<\mu <(n-2-\sigma)/2$, is an   entire solution  of inequality
(4)   such that $u(x)>v(x)$ a.e. in ${\mathbb R}^n$.
\end{example}

Combining Examples 1 and 2  gives the sharpness of Theorem \ref{thm1} in the
case when $1>q>0$.

In what follows, we consider  the case when $q\geq 1$.

\begin{theorem}\label{thm3}
Let $n\geq2$, $q>1$, $\nu \in (0,1)\cap (0,q-1)$,
$p=2(q-\nu)/(q-1)$, let $L$ be a differential operator defined by
\eqref{(1)} in $ {\mathbb R}^n$, let
\begin{eqnarray}\label{(21)}
\liminf\limits_{R\to\infty}{\mathrm {cap}}_{L,p} (\overline
{\mathcal B}_{R/2},{\mathbb R}^n\setminus {\mathcal B}_R; {\mathbb
R}^n )<\infty,
\end{eqnarray} and  $(u, v)$ be an entire solution of 
inequality \eqref{(4)} in ${\mathbb R}^n$ such that $u\geq v$ a.e. in ${\mathbb R}^n$. Then
$u = v$ a.e. in ${\mathbb R}^n$.
\end{theorem}

Now,  we give  conditions on the coefficients of the operator $L$
and on the parameter $p$  which guarantee  that inequality \eqref{(21)}
holds.

\begin{proposition}\label{prop4}
Let $n\geq2$,  and let $L$ be a differential operator defined by  \eqref{(1)} in $ {\mathbb R}^n$ such that its coefficients
satisfy condition \eqref{(14)} for all sufficiently large $R$ with some constants
$A>0$ and $n-2>\sigma>-2$. Then inequality \eqref{(21)}
holds for any $p\geq {2n}/(\sigma+2)$.
\end{proposition}

\begin{proposition}\label{prop5}
Let $n\geq2$,   and let $L$ be a differential operator defined by
\eqref{(1)} in $ {\mathbb R}^n$ such that its coefficients
 satisfy condition \eqref{(14)} for all sufficiently large $R$ with  some constants
$A>0$ and $n-2>\sigma> -2$. Assume also that $q>1$, $\nu
\in (0,1)\cap (0,q-1)$ and $q\leq (n-\nu(\sigma+2))/(n-\sigma -2)$.
Then inequality \eqref{(21)} holds for $p=2(q-\nu)/(q-1)$.
\end{proposition}

The following example shows the sharpness of Propositions 4 and 5.

\begin{example}
Let $n\geq 2$ and $p>1$. Consider the differential operator $L$
defined by \eqref{(17)} in ${\mathbb R}^n$, with some
constant $\sigma$. As in Example 1, using the standard arguments in capacity
theory (see, e.g., [2, p.178] or [8, p.12]), it is not difficult to
verify  that the $(L,p)$-capacity, associated with the operator $L$,
of the  condenser $(\overline {\mathcal B}_{R/2},{\mathbb
R}^n\setminus {\mathcal B}_R; {\mathbb R}^n )$ satisfies the
two-sided inequality
\begin{eqnarray}\label{(22)}
C_4 R^{(2n-p(\sigma+2))/2}\geq
{\mathrm {cap}}_{L,p} (\overline {\mathcal B}_{R/2},{\mathbb
R}^n\setminus {\mathcal B}_R; {\mathbb R}^n)\geq C_3
R^{(2n-p(\sigma+2))/2}
\end{eqnarray}
for all  $R>1$, with the same constant $\sigma$  as in \eqref{(17)} and some positive constants $C_3$, $C_4$  which depend only on
$n$, $p$ and $\sigma$. Moreover, for any $\sigma$ such that
$n-2>\sigma>-2$, inequality \eqref{(22)} yields that  if $p\geq 2n/(\sigma+2)$, then  for all $R>1$, the
$(L,p)$-capacity of the condenser $({\overline {\mathcal
B}}_{R/2},{\mathbb R}^n\setminus {\mathcal B}_R; {\mathbb R}^n )$ is
bounded above by a constant which depends only on $n$, $p$  and
$\sigma$, whereas $({\overline {\mathcal
B}}_{R/2},{\mathbb R}^n\setminus {\mathcal B}_R; {\mathbb R}^n )$ tends to
infinity as $R\to \infty$ if $p<2n/(\sigma+2)$. Thus, indeed, Proposition \ref{prop4} is sharp.

Moreover, inequality (22), for any $q>1$ and   any $\nu \in (0,1)\cap
(0,q-1)$ with $p=2(q-\nu)/(q-1)$   as in Theorem \ref{thm3} and
Proposition \ref{prop5},  can be rewritten in the form
\begin{align}\label{(23)}
	C_4 R^\frac{(n-\sigma-2)(q-(n-\nu(\sigma+2))}{(n-\sigma-2))(q-1)}&\geq \nonumber \\
	{\mathrm {cap}}_{L,p} &({\overline {\mathcal B}}_{R/2},{\mathbb
	R}^n\setminus {\mathcal B}_R; {\mathbb R}^n)\geq C_3
	R^\frac{(n-\sigma-2)(q- (n-\nu(\sigma+2))}{({n-\sigma-2}))(q-1)}
\end{align}
for all  $R>1$, with the same constants $\sigma$,  $C_3$ and  $C_4$ 
as in inequality \eqref{(22)}. It is easy to verify  using inequality \eqref{(23)} that for all $R>1$  if $1<q\leq  
({n-\nu(\sigma+2)})/(n-\sigma-2)$, then  for
any $\sigma$ such that $n-2>\sigma>-2$, the $(L,p)$-capacity  of the condenser $({\overline {\mathcal
B}}_{R/2},{\mathbb R}^n\setminus {\mathcal B}_R; {\mathbb R}^n )$ is
bounded above by a constant which depends only on $n$, $p$ and
$\sigma$. Whereas $({\overline {\mathcal
B}}_{R/2},{\mathbb R}^n\setminus {\mathcal B}_R; {\mathbb R}^n )$  tends to infinity as $R\to
\infty$ if  $q\geq  n /(n-\sigma-2)$. Thus, indeed, Proposition \ref{prop5} is sharp.
\end{example}

Since the positive  parameter $\nu$ in Proposition \ref{prop5}  may be chosen
arbitrarily small, the next result follows directly from Theorem \ref{thm3}
and Proposition \ref{prop5}.

\begin{theorem}\label{thm4}
Let $n\geq 2$. Assume that $L$ is a differential operator defined by
 \eqref{(1)} in $ {\mathbb R}^n$ such that its coefficients
satisfy condition \eqref{(14)} for all sufficiently large $R$ with some constants
$A>0$ and $n-2>\sigma>-2$. Let $1<q< n/(n-\sigma
-2)$,   and  $(u, v)$ be an entire solution of inequality
\eqref{(4)} in ${\mathbb R}^n$ such that $u\geq v$ a.e. in ${\mathbb R}^n$. Then $u = v$ a.e. 
in ${\mathbb R}^n$.
\end{theorem}

\begin{remark}
The case $\sigma\geq n-2$ is covered by Theorem
	2.
\end{remark}
 
The  next example demonstrates the
sharpness of the hypothesis $\sigma>-2$ in
Theorem \ref{thm4}, as well as in Propositions 4--8 and Theorems 6 and 8
given  below.

\begin{example}
Let $n\geq 2$ and  $q>1$.  Consider the differential operator $L$
from Example 1 defined by  \eqref{(17)} with $\sigma\leq -2$.
It is easy to verify that the pair $(u,v)$ of the functions
\begin{eqnarray*}\label{(24)}
u(x)=\alpha (1+|x|^2)^{-\mu}\quad {\rm{and}}\quad v(x)=0,
\end{eqnarray*}
where  $(n-\sigma-2)/2>\mu> -(\sigma+2)/2$ and $\alpha$ is  a
sufficiently small positive constant which depends only on
$\mu$, is an entire solution of  inequality \eqref{(4)}  such that $u>v$ a.e. in ${\mathbb R}^n$.
\end{example}

To complete our study in the  case $q>1$, consider  
$1<q\leq n/(n-\sigma-2)$ and  introduce  the quantity
\begin{eqnarray*}\label{(25)}
{\frak C}_{L, p_1,p_2}(R):= \left({\mathrm {cap}}_{L,p_1}
({\overline {\mathcal B}}_R,{\mathbb R}^n\setminus {\mathcal
B}_{2R}; {\mathbb R}^n )\right)^{1/2} \left({\mathrm {cap}}_{L,p_2}
(\overline {\mathcal B}_{R/2},{\mathbb R}^n\setminus {\mathcal B}_R;
{\mathbb R}^n )\right)^{1/p_2},
\end{eqnarray*}
which includes both  the $(L,p_1)$-capacity of the condenser $\left
({\overline {\mathcal B}}_R,{\mathbb R}^n\setminus {\mathcal
B}_{2R}; {\mathbb R}^n )\right)$ and the $(L,p_2)$-capacity of the
condenser $\left (\overline {\mathcal B}_{R/2},{\mathbb
R}^n\setminus {\mathcal B}_R; {\mathbb R}^n )\right)$ for all $R>0$,
where $p_1>1$ and $p_2>1$.

\begin{theorem}\label{thm5}
Let $n\geq2$, $q>1$, $\nu \in (0,1)\cap (0,q-1)$,
$p_1=2(q-\nu)/(q-1)$,  $p_2=2q/(q-1-\nu)$. Assume that $L$ is a differential
operator defined by \eqref{(1)} in $ {\mathbb R}^n$. Let
\begin{eqnarray}\label{(26)}
\liminf\limits_{R\to\infty}{\frak C}_{L, p_1,p_2}(R)<\infty,
\end{eqnarray}
and $(u, v)$ be an entire solution of inequality \eqref{(4)} in
${\mathbb R}^n$ such that $u\geq v$ a.e. in ${\mathbb R}^n$. Then $u = v$ in
${\mathbb R}^n$.
\end{theorem}

Before  imposing  conditions  on the coefficients of the operator $L$
and on the parameters $p_1$ and $p_2$ which would provide 
inequality \eqref{(26)}, we show by the following example that the
parameter $q$  must indeed satisfy the condition $1<q\leq n/(n-\sigma-2)$.

\begin{example} Let $n\geq 2$. Consider  the differential operator $L$ from Example 1 defined by  \eqref{(17)} with $n-2>\sigma>-2$, and $q> n/
(n-\sigma -2)$. It is easy to verify that  the pair $(u,v)$ of the
functions
\begin{eqnarray*}\label{(27)}
u(x)=\alpha (1+|x|^2)^{-\mu}\quad {\rm{and}}\quad v(x)=0,
\end{eqnarray*}
where $(n-2-\sigma)/2>\mu\geq (\sigma+2)/(2(q-1))$ and $\alpha$ is a
sufficiently small positive number which depends only on
$\mu$ is an entire  solution  of inequality \eqref{(4)} such that $u>v$ a.e. in $\mathbb{R}^n$.
\end{example}

\begin{proposition}\label{prop6}
Let $n\geq 2$ and  $L$ be a differential operator defined by
\eqref{(1)} in $ {\mathbb R}^n$  such that its coefficients
satisfy condition \eqref{(14)} for all sufficiently large $R$, with some constants
$A>0$, $n-2>\sigma> -2$. Furthermore, assume $\nu \in
(0,1)\cap (0,(\sigma+2)/ (n-\sigma -2))$. Then inequality \eqref{(26)}
holds for any $p_1\geq 2(n-\nu(n-\sigma-2))/(\sigma+2)$ and $p_2\geq
2n/(\sigma+2-\nu(n-\sigma-2))$.
\end{proposition}

Furthermore, we have
\begin{proposition}\label{prop7}
 Let $n\geq 2$  and $L$ be a differential operator defined
by \eqref{(1)} in $ {\mathbb R}^n$ such that its
coefficients satisfy condition \eqref{(14)} for all sufficiently large
$R$, with  some constants $A>0$, $\sigma$. Furthermore, assume that either
$n-2>\sigma> -2$,  $1<q\leq n/(n-\sigma-2)$ and $\nu \in (0,1)\cap
(0,q-1)$,  or   $\sigma\geq n-2$,  $q>1$ and $\nu \in (0,1)\cap
(0,q-1)$.  Then inequality \eqref{(26)} holds  for $p_1=2(q-\nu)/(q-1)$
and $p_2=2q/(q-1-\nu)$.
\end{proposition}

The following example shows  the sharpness of Propositions 6
and 7.

\begin{example}
Let $n\geq 2$ and  $p_2>p_1>2$. Consider the differential operator
$L$ defined by \eqref{(17)} in ${\mathbb R}^n$ with some constant $\sigma$. From \eqref{(22)}, it follows that the two-sided inequalities
\begin{eqnarray}\label{(28)}
C_6 R^{(2n-p_1(\sigma+2))/4}\geq (
{\mathrm {cap}}_{L,p_1} (\overline {\mathcal B}_{R/2},{\mathbb
R}^n\setminus {\mathcal B}_R; {\mathbb R}^n))^{1/2} \geq C_5
R^{(2n-p_1(\sigma+2))/4}
\end{eqnarray}
and
\begin{eqnarray}\label{(29)}
C_8 R^{(2n-p_2(\sigma+2))/(2p_2)}\geq
({\mathrm {cap}}_{L,p_2} ({\overline {\mathcal B}}_R,{\mathbb
R}^n\setminus {\mathcal B}_{2R}; {\mathbb R}^n))^{1/p_2}\geq C_7
R^{(2n-p_2(\sigma+2))/(2p_2)},
\end{eqnarray}
hold for all $R>1$, where $\sigma$ is the same constant as in \eqref{(17)}, whereas  $C_5$, $C_6$, $C_7$ and $C_8$ are some positive
constants which depend, possibly,  only on $n$, $p_1$, $p_2$ and
$\sigma$.

Furthermore, from  inequalities \eqref{(28)} and \eqref{(29)} it follows that  the two-sided
inequality
\begin{align}\label{(30)}
&C_{10} R^{(2n-p_1(\sigma+2))/4}
R^{(2n-p_2(\sigma+2))/(2p_2)}\geq \nonumber\\
&~~~~~{\frak C}_{L, p_1, p_2} (R)\geq  C_9
R^{(2n-p_1(\sigma+2))/4}R^{(2n-p_2(\sigma+2))/(2p_2)},
\end{align}
holds for all $R>1$, where $\sigma$ is the same constant as in  \eqref{(17)}, whereas $C_9$ and $C_{10}$ are some positive constants which
depend only on $n$, $p_1$, $p_2$ and $\sigma$.

From \eqref{(30)}, under the assumptions as in Proposition \ref{prop6}, that is for any $\sigma$ such that $n-2>\sigma >-2$ and  any $\nu
\in (0,1)\cap (0,(\sigma+2))/ (n-\sigma -2)$, choosing  $p_1$ and
$p_2$ such that $p_1\geq 2(n-\nu(n-\sigma-2))/(\sigma+2)$ and
$p_2\geq 2n/(\sigma+2-\nu(n-\sigma-2))$, we
have 
\begin{eqnarray*}\label{(31)}
C_{11}\geq {\frak C}_{L, p_1, p_2} (R),~~~\forall R>1.
\end{eqnarray*}
Here $C_{11}$ is some positive constant
that depends only on $n$, $q$, $p_1$, $p_2$, $\sigma$ and $\nu$. Thus, indeed, Proposition \ref{prop6} is sharp.

Let us show the sharpness of Proposition \ref{prop7}. From \eqref{(30)}, under the assumptions as in Theorem
5 and Proposition \ref{prop7}, that is 
for any $\sigma$, any $q>1$ and any $\nu \in (0,1)\cap (0,q-1)$,
letting $p_1={2(q-\nu)/(q-1)}$ and $p_2=2q/(q-1-\nu)$, we have the two-sided inequality
\begin{align}\label{(32)}
&C_{13} R^{(2q-1-\nu)(q(n-\sigma-2)-n)/(2q(q-1))}\geq \nonumber
\\  &~~~~~~~~~~~~~~{\frak C}_{L, p_1, p_2} (R)\geq   C_{12} R^{(2q-1-\nu)
(q(n-\sigma-2)-n)/(2q(q-1))},
\end{align}
which holds for all $R>1$, with the same constant $\sigma$  as in  \eqref{(17)} and with some positive constants $C_{12}$, $C_{13}$ which
depend on $n$, $q$, $\sigma$ and $\nu$.

Inequality \eqref{(32)} implies that if $1<q\leq n/(n-\sigma-2)$, then for all $R>1$ and for any $\sigma$ such that $n-2>\sigma>-2$,  ${\frak C}_{L, p_1,p_2}(R)$ is bounded above by a constant
which depends only on $n$, $q$, $\sigma$  and $\nu$. Furthermore,
${\frak C}_{L, p_1,p_2}(R)$ tends
to infinity as $R\to \infty$ if $q> n/(n-\sigma-2)$.

 Moreover, from \eqref{(32)} it follows that for any $\sigma$ such that $\sigma\geq
n-2$ and for all $R>1$, ${\frak C}_{L, p_1,p_2}(R)$ is bounded above
by a constant which depends only on $n$, $q$, $\sigma$  and $\nu$.
\end{example}

Combining Examples 5 and 6 shows  the sharpness of Theorem \ref{thm5}.

The next result follows directly from Theorem \ref{thm5} and Proposition \ref{prop7}.

\begin{theorem}\label{thm6}
Assume that $n\geq 2$ and  $L$ is a differential operator defined by
 \eqref{(1)} in $ {\mathbb R}^n$ such that its coefficients
satisfy condition \eqref{(14)} for all sufficiently large $R$, with some constants
$A>0$ and $n-2>\sigma >-2$.  Let $1<q\leq n/(n-\sigma
-2)$ and $(u, v)$ be an entire solution of inequality \eqref{(4)}
in ${\mathbb R}^n$ such that $u\geq v$ a.e. in ${\mathbb R}^n$. Then
$u = v$ a.e. in ${\mathbb R}^n$.
\end{theorem}

Note that the case when $\sigma\geq n-2$ is covered by
Theorem \ref{thm2}.

\begin{remark}
	Example 5 shows that  the hypothesis in
	Theorem \ref{thm6}, which requires  that $1<q\leq n/ (n-\sigma -2)$,  is
	sharp. Moreover, combining Example 3 with Examples 4 and  5 gives
	that the result in Theorem \ref{thm3} is sharp, except for the case when
	$q=n/ (n-\sigma -2)$ which is covered  by Theorem \ref{thm5}.
\end{remark}

Finally, we complete our study by considering the case when $n\geq
2$ and $q=1$.

\begin{theorem}\label{thm7}
Assume that $n\geq2$, $q=1$, and  $L$ is a differential operator defined by
\eqref{(1)} in $ {\mathbb R}^n$. Let
\begin{eqnarray}\label{(33)}
\liminf\limits_{R\to\infty}{\mathrm {cap}}_{L,2} (\overline
{\mathcal B}_{R/2},{\mathbb R}^n\setminus {\mathcal B}_R; {\mathbb
R}^n )R^{-n}=0,
\end{eqnarray} 
and $(u, v)$ be an entire solution of inequality \eqref{(4)} in ${\mathbb R}^n$ such that $u\geq v$ a.e. in ${\mathbb R}^n$. Then
$u = v$ a.e. in ${\mathbb R}^n$.
\end{theorem}

Now,  we give a  condition on the coefficients of the operator $L$
which guarantees that equality \eqref{(33)} holds.

\begin{proposition}\label{prop8}
Let $n\geq2$,  and $L$ be a differential operator defined by
\eqref{(1)} in $ {\mathbb R}^n$  such that its coefficients
satisfy condition \eqref{(14)} for  all sufficiently large $R$ with some  constants
$A>0$ and  $\sigma > -2$. Then equality \eqref{(33)}
holds.
\end{proposition}

The following result is a simple corollary of Theorem \ref{thm7} and
Proposition \ref{prop8}.

\begin{theorem}\label{thm8}
Assume that $n\geq2$, $q=1$, and  $L$ is a differential operator defined by
\eqref{(1)} in $ {\mathbb R}^n$  such that its coefficients
 satisfy condition \eqref{(14)} for  all sufficiently large $R$ with some  constants
$A>0$ and  $\sigma > -2$. Let $(u, v)$ be an
entire solution of inequality \eqref{(4)} in ${\mathbb R}^n$ such that $u\geq v$ a.e. in ${\mathbb R}^n$. Then
$u = v$ a.e. in ${\mathbb R}^n$.
\end{theorem}

Combining Example 1 with the following example gives  the sharpness
of Theorems \ref{thm7}, \ref{thm8} and Proposition \ref{prop8}.

\begin{example}
Let $n\geq 2$ and  $q=1$.  Consider the differential operator $L$
from Example 1 defined by \eqref{(17)} with $\sigma\leq -2$. It
is easy to verify that the pair $(u,v)$ of the functions
\begin{eqnarray*}\label{(34)}
u(x)=(1+|x|^2)^{-\mu}\quad {\rm{and}}\quad v(x)=0,
\end{eqnarray*}
where $\mu$ is a constant such that $1/(2n)\leq \mu\leq
{-(\sigma+2)/2} $ for $\sigma \leq -2 -1/n$ and $
{(n-\sigma-2-\sqrt{(n-\sigma-2)^2 -4})/4\leq }\mu\leq
{(n-\sigma-2+\sqrt{(n-\sigma-2)^2 -4})/4}$ for $\sigma > -2 -1/n$,
is an entire solution of inequality \eqref{(4)} with the operator $L$ in
${\mathbb R}^n$ such that $u(x)>v(x)$ in ${\mathbb R}^n$.
\end{example}

The results in Theorems \ref{thm1}, \ref{thm3}, \ref{thm5} and \ref{thm7} are new; they are also new in
the case when entire solutions of inequalities \eqref{(4)}, \eqref{(5)}--\eqref{(6)}
in ${\mathbb R}^n$ belong to the function space $W^{1, q}_{\mathrm
{loc}}({\mathbb R}^n)$. Theorems \ref{thm2}, \ref{thm4}, \ref{thm6} and \ref{thm8} were proved in \cite{[6]};
in the present paper, we show that these theorems can be obtained as
simple corollaries of Theorems \ref{thm1}, \ref{thm3}, \ref{thm5} and \ref{thm7}. The result in Theorem
\ref{thm5} and its  proof  generalize and correct the result in Theorem \ref{thm1} and
its proof obtained in \cite{[5]}. We would like also to note that the
results obtained in this paper were motivated by results established
in \cite{[1]} and \cite{[4]}.

\section {Proofs}

\vspace{5 mm} \noindent {\bf Proof of Theorem \ref{thm1}}. Let $n\geq 2$,
$q>0$,  and $L$ be a  differential operator defined by
\eqref{(1)},   and let  $(u,v)$ be an entire  solution of inequality \eqref{(4)}
in ${\mathbb R}^n$  such that $u(x)\geq v(x)$ a.e. in ${\mathbb R}^n$. Then,  from \eqref{(8)} we
have the inequality
\begin{equation} \label{(35)}
\int\limits_{{\mathbb
R}^n}\sum\limits_{i,j=1}^n a_{ij}(u-v)_{x_i}\zeta_{x_j} dx \geq
\int\limits_{{\mathbb R}^n}(|u|^{q-1}u - |v|^{q-1}v)\zeta dx,
\end{equation}
which holds for every non-negative function $\zeta\in
C^{\infty}({\mathbb R}^n)$  with compact support.

Set $w(x)=u(x)-v(x)$ and let  $R$ and $\varepsilon$  be
positive numbers,  and  $\varphi$ be a  function such that
$\varphi\in {\rm{Lip}}_{\rm{loc}}({\mathbb R}^n)$, $\varphi=1$ on
$\overline {\mathcal B}_{R/2}$, $\varphi=0$ outside ${\mathcal
B}_R$, and $1\geq \varphi \geq 0$ in ${{\mathbb R}^n}$. Without loss of generality, we may substitute  the test function $\zeta
(x)=(w(x)+\varepsilon )^{-1} \varphi^2 (x)$ in \eqref{(35)}. Then integrating by parts, we obtain the inequality
\begin{eqnarray}\label{(36)}
&&2\int\limits_{{\mathcal B}_R\setminus {\mathcal
B}_{R/2}}\sum\limits_{i,j=1}^n
a_{ij}w_{x_i}\varphi_{x_j}(w+\varepsilon )^{-1}\varphi dx \geq
\nonumber \\
&&\int\limits_{{\mathcal B}_R}\sum\limits_{i,j=1}^n
a_{ij}w_{x_i}w_{x_j}(w+\varepsilon )^{-2}\varphi^2 dx +
\int\limits_{{\mathcal B}_R} (|u|^{q-1}u -
|v|^{q-1}v)(w+\varepsilon) ^{-1}\varphi^2 dx.
\end{eqnarray}

Estimating  the integral  on the left side of \eqref{(36)} by H\"{o}lder's
inequality we obtain
\begin{eqnarray}\label{(37)}
2\left(\int\limits_{{\mathcal B}_R\setminus {\mathcal B}_{R/2}}
\sum\limits_{i,j=1}^n a_{ij}w_{x_i}w_{x_j}(w+\varepsilon)^{-2}
\varphi^2 dx\right)^{1/2}\times \nonumber \\  \left(
\int\limits_{{\mathcal B}_R\setminus {\mathcal B}_{R/2}}
\sum\limits_{i,j=1}^n a_{ij}\varphi_{x_i}\varphi_{x_j} dx
\right)^{1/2} \geq \int\limits_{{\mathcal B}_R}\sum\limits_{i,j=1}^n
a_{ij}w_{x_i}w_{x_j}(w+\varepsilon )^{-2}\varphi^2 dx+\nonumber \\
\int\limits_{{\mathcal B}_R} (|u|^{q-1}u - |v|^{q-1}v)
(w+\varepsilon)^{-1}\varphi^{2} dx.
\end{eqnarray}
Since  both terms   on the right side of \eqref{(37)}
are non-negative, we have  the inequalities
\begin{eqnarray*}\label{(38)}
2\left(\int\limits_{{\mathcal B}_R\setminus {\mathcal B}_{R/2}}
\sum\limits_{i,j=1}^n a_{ij}w_{x_i}w_{x_j}(w+\varepsilon)^{-2}
\varphi^2 dx\right)^{1/2}\times \nonumber \\  \left(
\int\limits_{{\mathcal B}_R\setminus {\mathcal B}_{R/2}}
\sum\limits_{i,j=1}^n a_{ij}\varphi_{x_i}\varphi_{x_j} dx
\right)^{1/2} \geq \int\limits_{{\mathcal B}_R} (|u|^{q-1}u -
|v|^{q-1}v) (w+\varepsilon)^{-2}\varphi^2 dx
\end{eqnarray*}
and
\begin{eqnarray*} \label{(39)}
4 \int\limits_{{\mathcal B}_R\setminus {\mathcal B}_{R/2}} \sum\limits_{i,j=1}^n
a_{ij}\varphi_{x_i}\varphi_{x_j}dx \geq  \int\limits_{{\mathcal
B}_R} \sum\limits_{i,j=1}^n a_{ij}w_{x_i}w_{x_j}
(w+\varepsilon)^{-2}\varphi^2 dx,
\end{eqnarray*}
which then   yield the inequalities
\begin{eqnarray}\label{(40)}
2\left(\int\limits_{{\mathcal B}_R\setminus {\mathcal B}_{R/2}}
\sum\limits_{i,j=1}^n
a_{ij}w_{x_i}w_{x_j}(w+\varepsilon)^{-2}dx\right)^{1/2}\times
\nonumber \\  \left( \int\limits_{{\mathcal B}_R}
\sum\limits_{i,j=1}^n a_{ij}\varphi_{x_i}\varphi_{x_j} dx
\right)^{1/2} \geq \int\limits_{{\mathcal B}_{R/2}} (|u|^{q-1}u -
|v|^{q-1}v) (w+\varepsilon)^{-2}dx
\end{eqnarray}
and
\begin{eqnarray} \label{(41)}4 \int\limits_{{\mathcal B}_R} \sum\limits_{i,j=1}^n
a_{ij}\varphi_{x_i}\varphi_{x_j}dx \geq  \int\limits_{{\mathcal
B}_{R/2}} \sum\limits_{i,j=1}^n a_{ij}w_{x_i}w_{x_j}
(w+\varepsilon)^{-2}dx.
\end{eqnarray}

Minimizing  the left sides of \eqref{(40)} and \eqref{(41)}
over all functions $\varphi (x)$ admissible in the definition of the
$(L,2)$-capacity of the condenser $(\overline {\mathcal B}_{R/2},
{\mathbb R}^n \setminus {\mathcal B}_R; {\mathbb R}^n)$, we obtain  the inequalities
\begin{eqnarray}\label{(42)}
2
\left(\int\limits_{{\mathcal B}_R\setminus {\mathcal B}_{R/2}}
\sum\limits_{i,j=1}^n
a_{ij}w_{x_i}w_{x_j}(w+\varepsilon)^{-2}dx\right)^{1/2}\times
\nonumber\\
\left( {\mathrm {cap}}_{L,2} (\overline{\mathcal B}_{R/2}, {\mathbb
R}^n \setminus {\mathcal B}_R; {\mathbb R}^n )\right)^{1/2} \geq
\int\limits_{{\mathcal B}_{R/2}} (|u|^{q-1}u - |v|^{q-1}v)
(w+\varepsilon)^{-2}dx
\end{eqnarray}
and
\begin{eqnarray}\label{(43)}
4{\mathrm {cap}}_{L,2} (\overline {\mathcal B}_{R/2}, {\mathbb R}^n
\setminus {\mathcal B}_R; {\mathbb R}^n )\geq \int\limits_{{\mathcal
B}_{R/2}} \sum\limits_{i,j=1}^n a_{ij}w_{x_i}w_{x_j}
(w+\varepsilon)^{-2} dx.
\end{eqnarray}
Moreover, since  by one of the  hypotheses of Theorem \ref{thm1} there exists
a non-negative number $\Gamma$ and  an increasing sequence of
positive numbers $R_k$ such that ${R_k\to \infty}$ and
\begin{eqnarray} \label{(44)}
{\mathrm {cap}}_{L,2} ({\overline {\mathcal B}}_{R_k/2}, {\mathbb
R}^n \setminus {\mathcal B}_{R_k}; {\mathbb R}^n )\to \Gamma
\end{eqnarray}
as $R_k\to \infty$,  then, from  \eqref{(43)} and \eqref{(44)},  we have the
inequality
\begin{eqnarray}\label{(45)}
 \int\limits_{{\mathcal B}_{R_k/2}}
\sum\limits_{i,j=1}^n a_{ij}w_{x_i}w_{x_j} (w+\varepsilon)^{-2} dx
\leq 4 \Gamma,
\end{eqnarray}
which holds as $R_k\to \infty$. Due to 
condition \eqref{(3)}, the quantity
\begin{eqnarray*}\label{(46)}
 H(R):= \int\limits_{{\mathcal B}_R}
\sum\limits_{i,j=1}^n a_{ij}w_{x_i}w_{x_j} (w+\varepsilon)^{-2} dx
\end{eqnarray*}
increases  monotonically with respect to $R$. Hence by inequality
\eqref{(45}, which holds as $R_k\to \infty$, we derive the inequality
\begin{eqnarray*}\label{(47)}
 \int\limits_{{\mathcal B}_R}
\sum\limits_{i,j=1}^n a_{ij}w_{x_i}w_{x_j} (w+\varepsilon)^{-2}
dx\leq 4 \Gamma,
\end{eqnarray*}
which holds for all $R>0$. Due to the monotonicity of
$H(R)$ with respect to $R$, $H(R)$ has a limit as
$R\to\infty$, bounded by the constant $4\Gamma$, 
namely,
\begin{eqnarray}\label{(48)}
\lim\limits_{R\to \infty}
\int\limits_{{\mathcal B}_R} \sum\limits_{i,j=1}^n
a_{ij}w_{x_i}w_{x_j} (w+\varepsilon)^{-2} dx \leq 4 \Gamma.
\end{eqnarray}
From \eqref{(48)},   again due to the monotonicity of $H(R)$ with respect to
$R$, we obtain the equality
\begin{eqnarray*}\label{(49)}
\lim\limits_{R\to \infty}
 \int\limits_{{\mathcal B}_R\setminus {\mathcal B}_{R/2}}
\sum\limits_{i,j=1}^n a_{ij}w_{x_i}w_{x_j} (w+\varepsilon)^{-2} = 0,
\end{eqnarray*}
and, in particular,
\begin{eqnarray}\label{(50)}
\lim\limits_{R_k\to \infty} \int\limits_{{\mathcal B}_{R_k}\setminus
{\mathcal B}_{R_k/2}} \sum\limits_{i,j=1}^n a_{ij}w_{x_i}w_{x_j}
(w+\varepsilon)^{-2} = 0,
\end{eqnarray}
where $R_k$ is the same sequence as in \eqref{(44)}.

Since $|u|^{q-1}u\geq  |v|^{q-1}v$ in ${\mathbb R}^n$, we
derive  that the quantity
\begin{eqnarray*}\label{(51)}
I(R):= \int\limits_{{\mathcal B}_{R/2}} (|u|^{q-1}u - |v|^{q-1}v)
(w+\varepsilon)^{-2}dx,
\end{eqnarray*}
which is the right side of \eqref{(42)}, increases monotonically with
respect to $R$ and, therefore,  has a limit as $R\to \infty$, which,
generally speaking, can be equal to infinity.

Now, taking in \eqref{(42)}  the same sequence $R=R_k$  as in  \eqref{(44)}  and
then passing to the limit as $R_k\to \infty$, we obtain, due to \eqref{(44)}
and \eqref{(50)}, the equality
\begin{eqnarray}\label{(52)}
\lim\limits_{R_k\to \infty}\int\limits_{{\mathcal B}_{R_k/2}}
(|u|^{q-1}u - |v|^{q-1}v) (w+\varepsilon)^{-2}dx=0.
\end{eqnarray}

Finally, from \eqref{(52)}, since $|u|^{q-1}u\geq |v|^{q-1}v$
in ${\mathbb R}^n$, we deduce that $u(x)=v(x)$ almost everywhere in
${\mathbb R}^n$, and this concludes the proof of Theorem \ref{thm1}.

\vspace{5 mm} \noindent \textbf {Proof of Theorem \ref{thm3}}. Let $n\geq2$,
$q>1$, let $L$ be a differential operator defined by
\eqref{(1)}, and let $(u ,v)$ be an entire solution of the inequality \eqref{(4)} in
${\mathbb R}^n$ such that $u(x)\geq v(x)$. Using the algebraic
inequality
\begin{equation}\label{(53)}(|u|^{q-1}u-|v|^{q-1}v)(u-v)\geq c_1 |u-v|^{q+1},
\end{equation}
which, for  every $q\geq 1$ and a some positive constant $c_1$ 
depending  only on $q$, holds for all real numbers $u$ and $v$, we
obtain from \eqref{(8)} the inequality
\begin{equation}\label{(54)}
\int\limits_{{\mathbb R}^n} \sum\limits_{i,j=1}^n a_{ij}
(u-v)_{x_i}\zeta_{x_j} dx \geq c_1 \int\limits_{{\mathbb
R}^n}(u-v)^q\zeta dx,
\end{equation}
which holds for every non-negative function $\zeta\in
C^{\infty}({\mathbb R}^n)$ with compact support.

Set $w(x)=u(x)-v(x)$, and let $R$ and $\varepsilon$  be
positive numbers, $\nu$ is a positive number such that $\nu \in
(0, 1)\cap (0,q-1)$. Let $\varphi$ be a function such that
$\varphi\in {\rm{Lip}}_{\rm{loc}}({\mathbb R}^n)$, $\varphi=1$ on
$\overline {\mathcal B}_{R/2}$,  $\varphi=0$ outside ${\mathcal
B}_R$, and $1\geq \varphi \geq 0$ in ${{\mathbb R}^n}$. 
Without loss of generality, we may substitute  the test function  $\zeta
(x)=(w(x)+\varepsilon)^{-\nu}\varphi^s (x)$ in \eqref{(54)}. Hence, we get the inequality
\begin{eqnarray}\label{(55)}
&&s\int\limits_{{\mathcal B}_R\setminus {\mathcal
B}_{R/2}}\sum\limits_{i,j=1}^n a_{ij}
w_{x_i}\varphi_{x_j}(w+\varepsilon )^{-\nu}\varphi^{s-1} dx \geq
\nonumber \\
&&\nu \int\limits_{{\mathcal B}_R}\sum\limits_{i,j=1}^n
a_{ij}w_{x_i}w_{x_j} (w+\varepsilon )^{-\nu -1}\varphi^s dx +
c_1\int\limits_{{\mathcal B}_R} w^q (w+\varepsilon) ^{-\nu}\varphi^s
dx.
\end{eqnarray}

Estimating the integrand on the left side of \eqref{(55)} by Cauchy's
inequality, we obtain  the inequality
\begin{eqnarray}\label{(56)}
\int\limits_{{\mathcal B}_R\setminus {\mathcal B}_{R/2}} s \left (
\sum\limits_{i,j=1}^n a_{ij}w_{x_i}w_{x_j}\right)^{1/2} \left (
\sum\limits_{i,j=1}^n
a_{ij}\varphi_{x_i}\varphi_{x_j}\right)^{1/2} (w+\varepsilon)^{-\nu}\varphi^{s-1} dx \geq \nonumber \\
\nu \int\limits_{{\mathcal B}_R}\sum\limits_{i,j=1}^n
a_{ij}w_{x_i}w_{x_j} (w+\varepsilon )^{-\nu -1}\varphi^s dx +
c_1\int\limits_{{\mathcal B}_R} w^q (w+\varepsilon) ^{-\nu}\varphi^s
dx.
\end{eqnarray}
Further estimating the integrand on the left side of \eqref{(56)} by
Young's inequality, we arrive at  the inequality
\begin{eqnarray}\label{(57)}
\frac \nu 2  \int\limits_{{\mathcal B}_R\setminus {\mathcal
B}_{R/2}}\sum\limits_{i,j=1}^n a_{ij}w_{x_i}w_{x_j}
(w+\varepsilon)^{-\nu-1}\varphi^s dx+ \nonumber
\\ c_2\int\limits_{{\mathcal B}_R\setminus {\mathcal B}_{R/2}} \sum\limits_{i,j=1}^n
a_{ij}\varphi_{x_i}\varphi_{x_j}
(w+\varepsilon)^{1-\nu}\varphi^{s-2}
dx\geq \nonumber \\
\nu \int\limits_{{\mathcal B}_R}\sum\limits_{i,j=1}^n
a_{ij}w_{x_i}w_{x_j} (w+\varepsilon )^{-\nu -1}\varphi^s dx +
c_1\int\limits_{{\mathcal B}_R} w^q (w+\varepsilon) ^{-\nu}\varphi^s
dx.
\end{eqnarray}
Here and in what follows in the proof of Theorem \ref{thm3},  we use the
symbols $c_i, i=2, \dots,$ to denote positive  constants depending
possibly on $n$, $q$ or $\nu$ but not on $R$ or  $\varepsilon$.

Due to condition \eqref{(3)}, from \eqref{(57)} we obtain the inequality
\begin{eqnarray}\label{(58)}
c_3\int\limits_{{\mathcal B}_R\setminus {\mathcal B}_{R/2}}
\sum\limits_{i,j=1}^n a_{ij}\varphi_{x_i}\varphi_{x_j}
(w+\varepsilon)^{1-\nu}\varphi^{s-2}dx \geq \int\limits_{{\mathcal
B}_R}w^q (w+\varepsilon)^{-\nu}\varphi^s dx.
\end{eqnarray}

Estimating the left side of \eqref{(58)} by  H\"{o}lder's inequality, we
arrive at  the inequality
\begin{eqnarray}\label{(59)}
c_3 \left(\int\limits_{{\mathcal B}_R\setminus {\mathcal
B}_{R/2}}\left(\sum\limits_{i,j=1}^n
a_{ij}\varphi_{x_i}\varphi_{x_j}\right)^{(q-\nu)/(q-1)}\varphi^{s-
2(q-\nu)/(q-1)} dx\right)^{(q-1)/(q-\nu)}\times\nonumber  \\
\left(\int\limits_{{\mathcal B}_R \setminus {\mathcal
B}_{R/2}}(w+\varepsilon)^{q-\nu}\varphi^s dx
\right)^{(1-\nu)/(q-\nu)} \geq  \int\limits_{{\mathcal
B}_R}w^q(w+\varepsilon)^{-\nu}\varphi^s dx.
\end{eqnarray}

Observe that $s=2(q-\nu)/(q-1)$. Hence setting
$p=2(q-\nu)/(q-1)$, we obtain   the inequality
\begin{eqnarray}\label{(60)}
&&c_3 \left(\int\limits_{{\mathcal B}_R\setminus {\mathcal
B}_{R/2}}\left(\sum\limits_{i,j=1}^n
a_{ij}\varphi_{x_i}\varphi_{x_j}\right)^{p/2}
dx\right)^{2/p}\times\nonumber  \\
&&\left(\int\limits_{{\mathcal B}_R \setminus {\mathcal
B}_{R/2}}(w+\varepsilon)^{p(q-1)/2}\varphi^s dx \right)^{(p-2)/p}
\geq \int\limits_{{\mathcal
B}_R}w^{p(q-1)/2+\nu}(w+\varepsilon)^{-\nu} \varphi^s dx.
\end{eqnarray}
Passing to the limit in \eqref{(60)} as $\varepsilon\to 0$, we obtain by
Lebesgue's theorem (see, e.g., \cite{[3]}, p.303) the
inequality
\begin{eqnarray*}\label{(61)}
&&c_3 \left(\int\limits_{{\mathcal B}_R\setminus {\mathcal
B}_{R/2}}\left(\sum\limits_{i,j=1}^n
a_{ij}\varphi_{x_i}\varphi_{x_j}\right)^{p/2}
dx\right)^{2/p}\times\nonumber  \\
&&\left(\int\limits_{{\mathcal B}_R \setminus {\mathcal
B}_{R/2}}w^{p(q-1)/2}\varphi^s dx \right)^{(p-2)/p} \geq
\int\limits_{{\mathcal B}_R}w^{p(q-1)/2}\varphi^s dx,
\end{eqnarray*}
which  yields 
\begin{eqnarray*}\label{(62)}
&&c_3 \left(\int\limits_{{\mathcal B}_R}\left(\sum\limits_{i,j=1}^n
a_{ij}\varphi_{x_i}\varphi_{x_j}\right)^{p/2}
dx\right)^{2/p}\times\nonumber  \\
&&\left(\int\limits_{{\mathcal B}_R \setminus {\mathcal
B}_{R/2}}w^{p(q-1)/2}dx \right)^{(p-2)/p} \geq
\int\limits_{{\mathcal B}_{R/2}}w^{p(q-1)/2}dx,
\end{eqnarray*}
and
\begin{eqnarray*}\label{(63)}
c_4 \int\limits_{{\mathcal B}_R}\left(\sum\limits_{i,j=1}^n
a_{ij}\varphi_{x_i}\varphi_{x_j}\right)^{p/2} dx \geq
\int\limits_{{\mathcal B}_{R/2}}w^{p(q-1)/2}dx.
\end{eqnarray*}
Minimizing  the left sides of these inequalities
over all functions $\varphi (x)$ admissible in the definition of the
$(L,p)$-capacity of the condenser $(\overline {\mathcal B}_{R/2},
{\mathbb R}^n \setminus {\mathcal B}_R; {\mathbb R}^n)$, we obtain
the inequalities
\begin{eqnarray}\label{(64)}
&&c_3 \left(\left( {\mathrm {cap}}_{L,p} ({\overline {\mathcal
B}}_{R/2}, {\mathbb R}^n \setminus {\mathcal B}_R; {\mathbb
R}^n)\right)^{p/2}\right)^{2/p}
\times\nonumber \\
 &&~~~~~~~\left(\int\limits_{{\mathcal B}_R \setminus
{\mathcal B}_{R/2}}w^{p(q-1)/2} dx \right)^{(p-2)/p} \geq
\int\limits_{{\mathcal B}_{R/2}}w^{p(q-1)/2}dx
\end{eqnarray}
and
\begin{eqnarray}\label{(65)}
c_4{\mathrm {cap}}_{L,p} (\overline {\mathcal B}_{R/2}, {\mathbb
R}^n \setminus {\mathcal B}_R; {\mathbb R}^n ) \geq
\int\limits_{{\mathcal B}_{R/2}}w^{ {p(q-1)}/2} dx.
\end{eqnarray}
Since by the  hypotheses of Theorem \ref{thm3} there exists a
non-negative number $\Lambda$ and  an increasing sequence of
positive numbers $R_k$ such that ${R_k\to \infty}$ and
\begin{eqnarray} \label{(66)}
{\mathrm {cap}}_{L,p} ({\overline {\mathcal B}}_{R_k/2}, {\mathbb
R}^n \setminus {\mathcal B}_{R_k}; {\mathbb R}^n )\to \Lambda
\end{eqnarray}
as $R_k\to \infty$, we obtain from \eqref{(65)} the inequality
\begin{eqnarray}\label{(67)}
\int\limits_{{\mathcal B}_{R_k/2}}w^{{p(q-1)}/2}dx \leq c_4\Lambda,
\end{eqnarray}
which holds as $R_k\to \infty$.  Moreover, the
quantity
\begin{eqnarray*} \label{(68)}
J(R):=
\int\limits_{{\mathcal B}_R}w^{{p(q-1)}/2}dx
\end{eqnarray*}
increases monotonically  with respect to $R$. Hence, from inequality
\eqref{(67)}, which holds as $R_k\to \infty$, we arrive at the inequality
\begin{eqnarray*}\label{(69)}
\int\limits_{{\mathcal B}_R}w^{{p(q-1)}/2}dx \leq c_4\Lambda,
\end{eqnarray*}
which holds for all $R>0$ and which, due to the monotonicity of
$J(R)$ with respect to $R$, yields that $J(R)$ has a limit, as $R\to
\infty$, bounded by the constant $c_4\Lambda$, namely,
\begin{eqnarray}\label{(70)}
\lim\limits_{R\to \infty}\int\limits_{{\mathcal
B}_R}w^{{p(q-1)}/2}dx \leq c_4\Lambda.
\end{eqnarray}
Hence,  due to the monotonicity of $J(R)$ with respect to $R$,
we have the equality
\begin{eqnarray*}\label{(71)}
\lim\limits_{R\to \infty} \int\limits_{{\mathcal B}_R\setminus
{\mathcal B}_{R/2}} w^{p(q-1)/2}dx = 0,
\end{eqnarray*}
and,  in paricular, the equality
\begin{eqnarray}\label{(72)}
\lim\limits_{R_k\to \infty} \int\limits_{{\mathcal B}_{R_k}\setminus
{\mathcal B}_{R_k/2}} w^{p(q-1)/2}dx = 0,
\end{eqnarray}
where $R_k$ is the same sequence as in \eqref{(66)}.

Moreover, observe that the right side of 
inequality \eqref{(64)}, which is  equal to $J(R/2)$,  increases monotonically
with respect to $R$ and, due to \eqref{(70)}, has a limit, as $R\to \infty$,
bounded from above  by $c_4\Lambda$. Let $R=R_k$ be the
same sequence as in \eqref{(66)}. Then  passing to the limit as $R_k\to
\infty$, due to \eqref{(66)}  and \eqref{(72)}, we obtain
\begin{eqnarray*}\label{(73)}
\lim\limits_{R_k\to \infty}\int\limits_{{\mathcal B}_{R_k/2}}
w^{p(q-1)/2} dx=0.
\end{eqnarray*}
This implies that $w(x)=0$ almost everywhere in
${{\mathbb R}^n}$, and thus $u(x)= v(x)$ almost
everywhere in ${{\mathbb R}^n}$. This concludes the proof of
Theorem \ref{thm3}.


\vspace{5 mm} \noindent \textbf {Proof of Theorem \ref{thm5}}. Let $n\geq2$,
$q>1$, $L$ be a differential operator defined by 
\eqref{(1)}, and $(u ,v)$ be an entire solution of the inequality \eqref{(4)} in
${\mathbb R}^n$ such that $u(x)\geq v(x)$. Then, as we show in the
proof of  Theorem \ref{thm3}, $(u,v)$  satisfies the inequality \eqref{(54)}. Let $w(x)=u(x)-v(x)$, and $R$ be a positive number, 
$\varphi$ be a function such that $\varphi\in
{\rm{Lip}}_{\rm{loc}}({\mathbb R}^n)$, $\varphi=1$ on $\overline
{\mathcal B}_{R/2}$,   $\varphi=0$ outside ${\mathcal B}_R$, and
$1\geq \varphi \geq 0$ in ${{\mathbb R}^n}$. Substituting 
the test function $\zeta(x)=\varphi^2(x)$ in \eqref{(54)} we obtain  the inequality
\begin{eqnarray*}\label{(74)}
2\int\limits_{{\mathbb R}^n} \sum\limits_{i,j=1}^n a_{ij}
w_{x_i}\varphi_{x_j} \varphi dx \geq c_1 \int\limits_{{\mathbb
R}^n}w^q\varphi^2 dx,
\end{eqnarray*}
where $c_1$ is the same constant as in \eqref{(53)}.

Let $\varepsilon$ be a   positive number and $\nu$ be a
positive number such that $\nu\in
 (0,1)\cap (q-1)$. Estimating the left side of \eqref{(73)} by H\"{o}lder's inequality, we
obtain the inequality
\begin{eqnarray}\label{(75)}
&&c_2\left(\int\limits_{{\mathcal B}_R\setminus {\mathcal
B}_{R/2}}\sum\limits_{i,j=1}^n
a_{ij}w_{x_i}w_{x_j} (w+\varepsilon)^{-1-\nu} dx\right)^{1/2}\times \nonumber \\
&&\left(\int\limits_{{\mathcal B}_R\setminus {\mathcal
B}_{R/2}}\sum\limits_{i,j=1}^n a_{ij}\varphi_{x_i}\varphi_{x_j}
(w+\varepsilon)^{1+\nu}\varphi^2 dx\right)^{1/2} \geq
\int\limits_{{\mathcal B}_R}w^q\varphi^2 dx.
\end{eqnarray}
Here and in what follows in the proof of Theorem \ref{thm5},  we use the
symbols $c_i, i=2, \dots,$ to denote positive  constants depending
possibly on $n$, $q$ or $\nu$ but not on $R$ or  $\varepsilon$.

Estimating the second term on the left side of \eqref{(75)} by
H\"{o}lder's inequality,  we arrive at  the inequality
\begin{eqnarray*}\label{(76)}
&&c_2\left(\int\limits_{{\mathcal B}_R\setminus {\mathcal
B}_{R/2}}\sum\limits_{i,j=1}^n
a_{ij}w_{x_i}w_{x_j} (w+\varepsilon)^{-1-\nu} dx\right)^{1/2}\times \nonumber \\
&&~~~~~~\left(\int\limits_{{\mathcal B}_R\setminus {\mathcal
B}_{R/2}}\left(\sum\limits_{i,j=1}^n
a_{ij}\varphi_{x_i}\varphi_{x_j}\right)^{q/(q-1-\nu)}\varphi^2
dx\right)^{(q-1-\nu)/(2q)}\times \nonumber \\
&&~~~~~~~~~~\left(\int\limits_{{\mathcal B}_R\setminus {\mathcal
B}_{R/2}}(w+\varepsilon)^{q}\varphi^{2} dx\right)^{(1+\nu)/(2q)}
\geq
 \int\limits_{{\mathcal B}_R}w^q\varphi^2 dx,
\end{eqnarray*}
which yields the inequality
\begin{eqnarray}\label{(77)}
&&c_2\left(\int\limits_{{\mathcal B}_R}\sum\limits_{i,j=1}^n
a_{ij}w_{x_i}w_{x_j} (w+\varepsilon)^{-1-\nu} dx\right)^{1/2}\times \nonumber \\
&&~~~~~\left(\int\limits_{{\mathcal B}_R}\left(\sum\limits_{i,j=1}^n
a_{ij}\varphi_{x_i}\varphi_{x_j}\right)^{q/(q-1-\nu)}
dx\right)^{(q-1-\nu)/(2q)}\times  \nonumber \\
&&~~~~~~~~\left(\int\limits_{{\mathcal B}_R\setminus {\mathcal
B}_{R/2}}(w+\varepsilon)^{q}\varphi^{2} dx\right)^{(1+\nu)/(2q)}
\geq \int\limits_{{\mathcal B}_R}w^q\varphi^2 dx.
\end{eqnarray}
Now, we estimate  the first  cofactor on the left side of \eqref{(77)}. To
this end, as above, we use the inequality \eqref{(54)}. Let $\psi$ be a
function such that ${\psi\in {\rm{Lip}}_{\rm{loc}}({\mathbb R}^n)}$,
$\psi=1$ on ${\overline {\mathcal B}}_R$,   $\psi=0$ outside
${\mathcal B}_{2R}$, and $1\geq \psi \geq 0$ in ${{\mathbb R}^n}$.
Without loss of generality, we may substitute  the test function  $\zeta
(x)=(w(x)+\varepsilon)^{-\nu}\psi^s (x)$ in \eqref{(54)}, where
$s= 2(q-\nu)/(q-1)$. Then integrating by parts, we get the inequality
\begin{eqnarray}\label{(78)}
&&s\int\limits_{{\mathcal B}_{2R}}\sum\limits_{i,j=1}^n a_{ij}
w_{x_i}\psi_{x_j}(w+\varepsilon)^{-\nu}\psi^{s-1} dx \geq
\nonumber \\
&&\nu \int\limits_{{\mathcal B}_{2R}}\sum\limits_{i,j=1}^n
a_{ij}w_{x_i}w_{x_j} (w+\varepsilon)^{-\nu -1}\psi^s dx +
c_1\int\limits_{{\mathcal B}_{2R}} w^q (w+\varepsilon) ^{-\nu}\psi^s
dx.
\end{eqnarray}
Estimating the integrand on the left side of \eqref{(78)} by  Cauchy's
inequality we obtain  the inequality
\begin{eqnarray}\label{(79)}
&&\int\limits_{{\mathcal B}_{2R}} s \left ( \sum\limits_{i,j=1}^n
a_{ij}w_{x_i}w_{x_j}\right)^{1/2} \left ( \sum\limits_{i,j=1}^n
a_{ij}\psi_{x_i}\psi_{x_j}\right)^{1/2} (w+\varepsilon)^{-\nu}\psi^{s-1} dx \geq \nonumber \\
&&\nu \int\limits_{{\mathcal B}_{2R}}\sum\limits_{i,j=1}^n
a_{ij}w_{x_i}w_{x_j} (w+\varepsilon)^{-\nu -1}\psi^s dx +
c_1\int\limits_{{\mathcal B}_{2R}} w^q (w+\varepsilon) ^{-\nu}\psi^s
dx.
\end{eqnarray}
Next, estimating the integrand on the left side of \eqref{(79)}   by Young's
inequality, we arrive at  the inequality
\begin{eqnarray*}\label{(80)}
\frac \nu 2  \int\limits_{{\mathcal B}_{2R}}\sum\limits_{i,j=1}^n
a_{ij}w_{x_i}w_{x_j} (w+\varepsilon)^{-\nu-1}\psi^s dx+
c_3\int\limits_{{\mathcal B}_{2R}} \sum\limits_{i,j=1}^n
a_{ij}\psi_{x_i}\psi_{x_j} (w+\varepsilon)^{1-\nu}\psi^{s-2}
dx\geq \nonumber \\
\nu \int\limits_{{\mathcal B}_{2R}}\sum\limits_{i,j=1}^n
a_{ij}w_{x_i}w_{x_j} (w+\varepsilon)^{-\nu -1}\psi^s dx +
c_1\int\limits_{{\mathcal B}_{2R}} w^q (w+\varepsilon) ^{-\nu}\psi^s
dx,
\end{eqnarray*}
which yields   the inequality
\begin{eqnarray}\label{(81)}
c_3\int\limits_{{\mathcal B}_{2R}} \sum\limits_{i,j=1}^n
a_{ij}\psi_{x_i}\psi_{x_j} (w+\varepsilon)^{1-\nu}\psi^{s-2}dx
\geq \nonumber \\
\frac \nu 2  \int\limits_{{\mathcal B}_{2R}}\sum\limits_{i,j=1}^n
a_{ij}w_{x_i}w_{x_j} (w+\varepsilon)^{-\nu-1}\psi^s dx+
c_1\int\limits_{{\mathcal B}_{2R}}w^q (w+\varepsilon)^{-\nu}\psi^s
dx.
\end{eqnarray}
Estimating the left side of \eqref{(81)} by  Young's  inequality we obtain
the inequality
\begin{eqnarray*}\label{(82)}
c_1\int\limits_{{\mathcal B}_{2R}}(w+\varepsilon)^{q-\nu}\psi^s dx +
c_4\int\limits_{{\mathcal B}_{2R}}\left(\sum\limits_{i,j=1}^n
a_{ij}\psi_{x_i}\psi_{x_j}\right)^{(q-\nu)/(q-1)}
\psi^{s-2(q-\nu)/(q-1)}dx \geq \nonumber
\\ \frac \nu 2 \int\limits_{{\mathcal B}_{2R}}\sum\limits_{i,j=1}^n
a_{ij}w_{x_i}w_{x_j} (w+\varepsilon)^{-\nu-1}\psi^s dx+
c_1\int\limits_{{\mathcal B}_{2R}}w^q(w+\varepsilon)^{-\nu}\psi^s
dx,
\end{eqnarray*}
which, by $s=2(q-\nu)/(q-1)$,  yields the inequality
\begin{eqnarray}\label{(83)}
&&c_1 \int\limits_{{\mathcal B}_{2R}}(w+\varepsilon)^{q-\nu}\psi^s dx
- c_1 \int\limits_{{\mathcal B}_{2R}}w^q(w+\varepsilon)^{-\nu}\psi^s
dx + \nonumber
\\ &&c_4\int\limits_{{\mathcal B}_{2R}}\left(\sum\limits_{i,j=1}^n
a_{ij}\psi_{x_i}\psi_{x_j}\right)^{(q-\nu)/(q-1)} dx \geq
\frac{\nu}2 \int\limits_{{\mathcal B}_R}\sum\limits_{i,j=1}^n
a_{ij}w_{x_i}w_{x_j} (w+\varepsilon)^{-\nu-1}dx.
\end{eqnarray}
Now, we return to the inequality \eqref{(77)} and estimate the first
cofactor on the left side of this inequality by \eqref{(83)}. As a result,
we have the inequality
\begin{eqnarray}\label{(84)}
c_2 \left(2c_1{\nu}^{-1} \int\limits_{{\mathcal
B}_{2R}}(w+\varepsilon)^{q-\nu}\psi^s dx - 2c_1{\nu}^{-1}
\int\limits_{{\mathcal B}_{2R}}w^q(w+\varepsilon)^{-\nu}\psi^s dx +
\right. \nonumber \\
\left.  2c_4{\nu}^{-1} \int\limits_{{\mathcal
B}_{2R}}\left(\sum\limits_{i,j=1}^n
a_{ij}\psi_{x_i}\psi_{x_j}\right)^{(q-\nu)/(q-1)} dx
\right)^{1/2}\times \nonumber
\\
\left(\int\limits_{{\mathcal B}_R}\left(\sum\limits_{i,j=1}^n
a_{ij}\varphi_{x_i}\varphi_{x_j}\right)^{q/(q-1-\nu)}
dx\right)^{(q-1-\nu)/(2q)}\times \nonumber \\
\left(\int\limits_{{\mathcal B}_R\setminus {\mathcal
B}_{R/2}}(w+\varepsilon)^{q}\varphi^{2} dx\right)^{(1+\nu)/(2q)}
\geq \int\limits_{{\mathcal B}_R}w^q \varphi^{2} dx.
\end{eqnarray}
Passing to the limit in \eqref{(84)} as $\varepsilon\to 0$, 
 we obtain by Lebesgue's theorem (see, e.g., \cite{[3]}, p.303) the inequality
\begin{eqnarray*}\label{(85)}
&&c_5 \left(\int\limits_{{\mathcal B}_{2R}}\left(\sum\limits_{i,j=1}^n
a_{ij}\psi_{x_i}\psi_{x_j}\right)^{(q-\nu)/(q-1)} dx
\right)^{1/2}\times \nonumber
\\
&&~~~~~\left(\int\limits_{{\mathcal B}_R}\left(\sum\limits_{i,j=1}^n
a_{ij}\varphi_{x_i}\varphi_{x_j}\right)^{q/(q-1-\nu)}
dx\right)^{(q-1-\nu)/(2q)}\times \nonumber  \\
&&~~~~~~~~~~\left(\int\limits_{{\mathcal B}_R\setminus {\mathcal
B}_{R/2}}w^{q}\varphi^{2} dx\right)^{(1+\nu)/(2q)} \geq
\int\limits_{{\mathcal B}_{R}}w^{q}\varphi^{2} dx,
\end{eqnarray*}
which yields the inequalities
\begin{eqnarray}\label{(86)}
&&c_5 \left(\int\limits_{{\mathcal B}_{2R}}\left(\sum\limits_{i,j=1}^n
a_{ij}\psi_{x_i}\psi_{x_j}\right)^{(q-\nu)/(q-1)} dx
\right)^{1/2}\times \nonumber
\\
&&~~~~~~\left(\int\limits_{{\mathcal B}_R}\left(\sum\limits_{i,j=1}^n
a_{ij}\varphi_{x_i}\varphi_{x_j}\right)^{q/(q-1-\nu)}
dx\right)^{(q-1-\nu)/(2q)}\times \nonumber  \\
&&~~~~~~~~~~~~~~\left(\int\limits_{{\mathcal B}_R\setminus {\mathcal B}_{R/2}}w^{q}
dx\right)^{(1+\nu)/(2q)} \geq \int\limits_{{\mathcal
B}_{R/2}}w^{q}dx
\end{eqnarray}
and
\begin{eqnarray}\label{(87)}
&&c_5 \left(\int\limits_{{\mathcal B}_{2R}}\left(\sum\limits_{i,j=1}^n
a_{ij}\psi_{x_i}\psi_{x_j}\right)^{(q-\nu)/(q-1)} dx
\right)^{1/2}\times \nonumber
\\
&&~~\left(\int\limits_{{\mathcal B}_R}\left(\sum\limits_{i,j=1}^n
a_{ij}\varphi_{x_i}\varphi_{x_j}\right)^{q/(q-1-\nu)}
dx\right)^{(q-1-\nu)/(2q)}\geq  \left(\int\limits_{{\mathcal B}_{R/2}}w^{q}
dx\right)^{(2q-1-\nu)/(2q)}.
\end{eqnarray}

Taking $p_1=2(q-\nu)/(q-1)$,  $p_2=2q/(q-1-\nu)$ in \eqref{(86)} and \eqref{(87)}, and minimizing the left sides of these inequalities first over all
functions $\psi (x)$ admissible in the definition of the
$(L,p_1)$-capacity of the condenser $(\overline {\mathcal B}_{R},
{\mathbb R}^n \setminus {\mathcal B}_{2R}; {\mathbb R}^n)$ and then
over all functions $\varphi (x)$ admissible in the definition of the
$(L,p_2)$-capacity of the condenser $(\overline {\mathcal B}_{R/2},
{\mathbb R}^n \setminus {\mathcal B}_R; {\mathbb R}^n)$, we arrive
at the inequalities
\begin{eqnarray}\label{(88)}
c_5 \left({\mathrm {cap}}_{L,p_1} ({\overline {\mathcal B}}_{R},
{\mathbb R}^n \setminus {\mathcal B}_{2R}; {\mathbb R}^n)
\right)^{1/2} \left({\mathrm {cap}}_{L,p_2} ({\overline {\mathcal
B}}_{R/2},
{\mathbb R}^n \setminus {\mathcal B}_R; {\mathbb R}^n) \right)^{1/p_2}\times\nonumber\\
\left(\int\limits_{{\mathcal B}_R\setminus {\mathcal B}_{R/2}}w^{q}
dx\right)^{(1+\nu)/(2q)} \geq   \int\limits_{{\mathcal B}_{R/2}}w^q
dx,
\end{eqnarray}
and thus we get
\begin{eqnarray}\label{(89)}
 &&c_5 \left({\mathrm {cap}}_{L,p_1} ({\overline {\mathcal B}}_{R},
{\mathbb R}^n \setminus {\mathcal B}_{2R}; {\mathbb R}^n) \right)^{1/2}\times\nonumber\\
&&~~~~~~~~\left({\mathrm {cap}}_{L,p_2} ({\overline {\mathcal B}}_{R/2},
{\mathbb R}^n \setminus {\mathcal B}_R; {\mathbb R}^n)
\right)^{1/p_2} \geq \left(\int\limits_{{\mathcal B}_{R/2}}w^q
dx\right)^{(2q-1-\nu)/(2q)}.
\end{eqnarray}
Since by one of the  hypotheses of Theorem \ref{thm5} there exists a
non-negative number $\Upsilon$ and  an increasing sequence of
positive numbers $R_k$ such that $R_k\to \infty$ and
\begin{eqnarray}\label{(90)}
{\frak C}_{L, p_1, p_2} (R_k) \to \Upsilon
\end{eqnarray}
as $R_k\to \infty$, we have from \eqref{(89)}  the inequality
\begin{eqnarray}\label{(91)}
\int\limits_{{\mathcal B}_{R_k/2}}w^qdx \leq
c_6\Upsilon^{2q/(2q-1-\nu)},
\end{eqnarray}
which holds as $R_k\to \infty$. Observe that the quantity
\begin{eqnarray*} \label{(92)}
Q(R):=  \int\limits_{{\mathcal B}_R}w^qdx
\end{eqnarray*}
increases  monotonically  with respect to $R$. Hence, from inequality
\eqref{(91)}, which holds as $R_k\to \infty$, we derive the inequality
\begin{eqnarray*} \label{(93)} 
\int\limits_{{\mathcal B}_R}w^qdx\leq c_6\Upsilon^{2q/(2q-1-\nu)},
\end{eqnarray*}
which holds for all $R>0$. Due to the monotonicity of
$Q(R)$ with respect to $R$, this yields that $Q(R)$ has a limit, as $R\to
\infty$, bounded from above by the constant
$c_6\Upsilon^{2q/(2q-1-\nu)}$, namely,
\begin{eqnarray}\label{(94)}
\lim\limits_{R\to \infty} \int\limits_{{\mathcal B}_R}w^qdx \leq
c_6\Upsilon^{2q/(2q-1-\nu)}.
\end{eqnarray}
From \eqref{(94)}, again due to the monotonicity of $Q(R)$ with respect to
$R$, we have the equality
\begin{eqnarray}\label{(95)}
\lim\limits_{R\to \infty} \int\limits_{{\mathcal B}_R\setminus
{\mathcal B}_{R/2}}w^qdx =0,
\end{eqnarray}
and thus we have the equality
\begin{eqnarray*}\label{(96)}
\lim\limits_{R_k\to \infty} \int\limits_{{\mathcal B}_{R_k}\setminus
{\mathcal B}_{R_k/2}}w^qdx=0,
\end{eqnarray*}
where $R_k$ is the same sequence as in \eqref{(90)}.

Observe that the right side of \eqref{(88)}, which is equal to  $Q(R/2)$,  increases monotonically
with respect to $R$. Thus, due to \eqref{(94)}, it has a limit, as $R\to
\infty$, bounded from above  by $c_6\Upsilon^{2q/(2q-1-\nu)}$. Hence, taking  the same sequence $R=R_k$ as in \eqref{(90)} and passing to
the limit as $R_k\to \infty$, we obtain, due to \eqref{(90)}  and (95), the
equality
\begin{eqnarray}\label{(97)}
\lim\limits_{R_k\to \infty}\int\limits_{{\mathcal B}_{R_k/2}} w^q
dx=0.
\end{eqnarray}
Finally, from \eqref{(97)}, we deduce that $w(x)=0$  almost everywhere  in
${{\mathbb R}^n}$, which implies that $u(x)= v(x)$ almost
everywhere in ${{\mathbb R}^n}$, and this concludes the proof of
Theorem \ref{thm5}.

\vspace{5 mm} \noindent \textbf {Proof of Theorem \ref{thm7}}. We prove this
theorem by contradiction. Let $n\geq2$, $q=1$,   $L$ be a
differential operator defined by \eqref{(1)}, and  $(u ,v)$
be an entire solution of inequality \eqref{(4)} in ${\mathbb R}^n$ such
that $u(x)> v(x)$.  Then, as it has been shown in Theorem \ref{thm3}, $(u,v)$ satisfies inequality \eqref{(54)} with $q=1$ and $c_1=1$. Further, let $R$ and
$\varepsilon$ be positive numbers, and $\varphi$ be a function
such that $\varphi\in {\rm{Lip}}_{\rm{loc}}({\mathbb R}^n)$,
$\varphi=1$ on $\overline {\mathcal B}_{R/2}$,   $\varphi=0$ outside
${\mathcal B}_R$, and $1\geq \varphi \geq 0$ in ${{\mathbb R}^n}$.

Without loss of generality we may substitute  the function $\zeta
(x)=(w(x)+\varepsilon)^{-1}\varphi^2 (x)$ in \eqref{(54)}, where
$w(x)=u(x)-v(x)$. Then integrating by parts,  we have  the
inequality
\begin{eqnarray}\label{(98)}
&&2\int\limits_{{\mathcal B}_R\setminus {\mathcal
B}_{R/2}}\sum\limits_{i,j=1}^n a_{ij}
w_{x_i}\varphi_{x_j}(w+\varepsilon )^{-1}\varphi dx \geq
\nonumber \\
&&~~~~~~~~~~~~~\int\limits_{{\mathcal B}_R\setminus {\mathcal
B}_{R/2}}\sum\limits_{i,j=1}^n a_{ij}w_{x_i}w_{x_j} (w+\varepsilon
)^{-2}\varphi^2 dx + \int\limits_{{\mathcal B}_R} w (w+\varepsilon)
^{-1}\varphi^2 dx.
\end{eqnarray}
Estimating the integrand on the left side of \eqref{(98)} by Cauchy's
inequality, we obtain  the inequality
\begin{eqnarray}\label{(99)}
\int\limits_{{\mathcal B}_R\setminus {\mathcal B}_{R/2}} 2 \left (
\sum\limits_{i,j=1}^n a_{ij}w_{x_i}w_{x_j}\right)^{1/2} \left (
\sum\limits_{i,j=1}^n a_{ij}\varphi_{x_i}\varphi_{x_j}\right)^{1/2}
(w+\varepsilon)^{-1}\varphi  dx \geq \nonumber \\
 \int\limits_{{\mathcal B}_R}\sum\limits_{i,j=1}^n a_{ij}w_{x_i}w_{x_j}
(w+\varepsilon )^{-2}\varphi^2 dx + \int\limits_{{\mathcal B}_R} w
(w+\varepsilon) ^{-1}\varphi^2 dx.
\end{eqnarray}
Now estimating the integrand on the left side of \eqref{(99)} by
Young's inequality, we arrive at  the inequality
\begin{eqnarray}\label{(100)}
\frac 1 2  \int\limits_{{\mathcal B}_R\setminus {\mathcal
B}_{R/2}}\sum\limits_{i,j=1}^n a_{ij}w_{x_i}w_{x_j}
(w+\varepsilon)^{-2}\varphi^2 dx+ c_1\int\limits_{{\mathcal
B}_R\setminus {\mathcal B}_{R/2}} \sum\limits_{i,j=1}^n
a_{ij}\varphi_{x_i}\varphi_{x_j}
dx\geq \nonumber \\
\int\limits_{{\mathcal B}_R\setminus {\mathcal B}_{R/2}}
\sum\limits_{i,j=1}^n a_{ij}w_{x_i}w_{x_j} (w+\varepsilon
)^{-2}\varphi^2 dx + \int\limits_{{\mathcal B}_R} w (w+\varepsilon)
^{-1}\varphi^2 dx.
\end{eqnarray}
Here and in what follows in the proof of Theorem \ref{thm7},  we use the
symbols $c_i, i=1, \dots,$ to denote positive  constants depending
possibly on $n$ but not on $R$ or $\varepsilon$.

From \eqref{(100)} due to \eqref{(3)}  we have the
inequality
\begin{eqnarray*}\label{(101)}
c_2\int\limits_{{\mathcal B}_R\setminus {\mathcal B}_{R/2}}
\sum\limits_{i,j=1}^n a_{ij}\varphi_{x_i}\varphi_{x_j} dx \geq
\int\limits_{{\mathcal B}_R}w (w+\varepsilon)^{-1}\varphi^2 dx,
\end{eqnarray*}
which   yields the inequality
\begin{eqnarray}\label{(102)}
c_2\int\limits_{{\mathcal B}_R} \sum\limits_{i,j=1}^n
a_{ij}\varphi_{x_i}\varphi_{x_j} dx \geq \int\limits_{{\mathcal
B}_{R/2}}w (w+\varepsilon)^{-1} dx.
\end{eqnarray}
Passing to the limit in \eqref{(102)} as $\varepsilon\to 0$ by Lebesgue's theorem (see, e.g., \cite{[3]}, p.303) we derive 
\begin{eqnarray}\label{(103)}
c_2\int\limits_{{\mathcal B}_R} \sum\limits_{i,j=1}^n
a_{ij}\varphi_{x_i}\varphi_{x_j} dx \geq \int\limits_{{\mathcal
B}_{R/2}} dx.
\end{eqnarray}
Minimizing  the left side of \eqref{(103)}  over all
functions $\varphi (x)$ admissible in the definition of the
$(L,2)$-capacity of the condenser ${(\overline {\mathcal B}_{R/2},
{\mathbb R}^n \setminus {\mathcal B}_R; {\mathbb R}^n)}$, we obtain
the inequality
\begin{eqnarray}\label{(104)}
c_3{\mathrm {cap}}_{L,2} (\overline {\mathcal B}_{R/2}, {\mathbb
R}^n \setminus {\mathcal B}_R; {\mathbb R}^n) \geq R^n.
\end{eqnarray}
However, by a hypothesis of Theorem \ref{thm7} there exists an increasing
sequence of positive numbers $R_k\to \infty$ such that the equality
\begin{eqnarray} \label{(105)}
\lim\limits_{R_k\to \infty} {\mathrm {cap}}_{L,2}
({\overline {\mathcal B}}_{R_k/2}, {\mathbb R}^n \setminus {\mathcal
B}_{R_k}; {\mathbb R}^n )R_k^{-n}=0
\end{eqnarray} 
holds. This implies the desired contradiction, i.e. inequality \eqref{(104)} and 
equality \eqref{(105)} contradict  each other.


\vskip 5pt \noindent {\bf Proof of Proposition \ref{prop1}.} Let $n\geq 2$,
$p>1$,  and $L$ be a  differential operator defined by \eqref{(1)} in $ {\mathbb R}^n$. Using  the algebraic  inequality
\begin{eqnarray*}\label{(106)}
\sum\limits_{i,j=1}^n  a_{ij}(x){\xi}_i {\xi}_j \leq \left
(\sum\limits_{i,j=1}^n  a^2_{ij}(x) \right)^{\frac 12}
\left(\sum\limits_{i=1}^n {\xi}^2_i\right)^{\frac 14}
\left(\sum\limits_{j=1}^n {\xi}^2_j\right)^{\frac 14},
\end{eqnarray*} which
holds for all $\xi=(\xi_1, \dots, \xi_n)\in {\mathbb R}^n$ and
almost all $x\in {\mathbb R}^n$, we obtain  the inequality
\begin{eqnarray}\label{(107)}
\int\limits_{{\mathbb R}^n}\left(\sum\limits_{i,j=1}^n
a_{ij}(x){\varphi}_{x_i}{\varphi}_{x_j}\right)^{\frac p2}dx\leq
\sup_{x\in {\mathcal B}_R \setminus {\mathcal
B}_r}\left(\sum\limits_{i,j=1}^n a^2_{ij}(x)\right)^{\frac
p4}\int\limits_{{\mathbb R}^n}  |\nabla {\varphi}|^pdx,
\end{eqnarray}
which holds for all functions $\varphi (x)$ admissible in the
definition of the $(L,p)$-capacity of the condenser $(\overline
{\mathcal B}_{r},{\mathbb R}^n\setminus {\mathcal B}_R ;{\mathbb
R}^n)$, with $R>r>1$.

Minimizing first the left side and then the right side of
\eqref{(107)} over all functions $\varphi (x)$ admissible in the
definition of the $(L,p)$-capacity   of the condenser $(\overline
{\mathcal B}_{r},{\mathbb R}^n\setminus {\mathcal B}_R ;{\mathbb
R}^n)$, we obtain inequality \eqref{(12)} for any $R>r>1$.

\vskip 5pt \noindent {\bf Proof of Proposition \ref{prop2}.} Let $n\geq 2$,
$p>1$, and let $L$ be a  differential operator defined by the
relation \eqref{(1)} in $ {\mathbb R}^n$   and such that its coefficients
satisfy condition \eqref{(14)} for all sufficiently large $R$,  with
$A>0$ and $\sigma$ some constants. Then,  from \eqref{(12)} and \eqref{(14)}, we
obtain the inequality
\begin{eqnarray}\label{(108)}
{\mathrm {cap}}_{L,p} (\overline {\mathcal B}_{R/2},{\mathbb
R}^n\setminus {\mathcal B}_R ;{\mathbb R}^n)\leq A^{p/4} R^{-\sigma
p/2} {\mathrm {cap}}_p (\overline {\mathcal B}_{R/2},{\mathbb
R}^n\setminus {\mathcal B}_R ;{\mathbb R}^n),
\end{eqnarray} which holds for  all  sufficiently large $R$,
with  $A>0$ and  $\sigma$  the same  constants as in \eqref{(14)}. In turn,
from \eqref{(108)} and \eqref{(13)},  we derive  the inequality \eqref{(15)}, which holds
for all sufficiently large $R$, with $\sigma$ the same constant as
in \eqref{(14)} and $\hat C$ some positive constant which depends only on
$A$, $n$, $p$ and $\sigma$,   and  this concludes the proof of
Proposition 2.

\vskip 5pt \noindent {\bf Proof of Proposition \ref{prop3}.} Let $n\geq2$, and
$L$ be a differential operator defined by \eqref{(1)} in $
{\mathbb R}^n$  such that its coefficients  satisfy
condition \eqref{(14)} for  all sufficiently large $R$,  with  some  constants $A>0$ and $\sigma \geq n-2$. Then, in  \eqref{(15)}, letting $p=2$, we
derive the inequality
\begin{eqnarray}\label{(109)}
{\mathrm {cap}}_{L,2} (\overline {\mathcal B}_{R/2},{\mathbb
R}^n\setminus {\mathcal B}_R ;{\mathbb R}^n)\leq {\hat
C}R^{n-\sigma-2},
\end{eqnarray}
which holds for all sufficiently large $R$, with some
positive constant $\hat C$  which depends only on $A$, $n$ and $\sigma$. From
\eqref{(109)}, we obtain that for any ${\sigma \geq n-2}$ the $(L,
2)$-capacity of the condenser $(\overline {\mathcal
B}_{R/2},{\mathbb R}^n\setminus {\mathcal B}_R ;{\mathbb R}^n)$ is
bounded above by a constant  which depends only on $A$, $n$ and
$\sigma$, for all sufficiently large $R$. This concludes the
proof of Proposition \ref{prop3}.

\vskip 5pt \noindent {\bf Proof of  Proposition \ref{prop4}.} Let $n\geq2$,
and $L$ be a differential operator defined by \eqref{(1)}
in $ {\mathbb R}^n$  such that its coefficients satisfy 
condition \eqref{(14)} for all sufficiently large $R$,  with some constants $A>0$ and
$n-2>\sigma>-2$. Then, from \eqref{(15)} we
obtain that for any $p\geq 2n/(\sigma+2)$,   the $(L, p)$-capacity of
the condenser $({\overline {\mathcal B}}_{R/2},{\mathbb
R}^n\setminus {\mathcal B}_R ;{\mathbb R}^n)$ is bounded above by a
constant,  which depends only on $A$, $n$, $p$ and $\sigma$, for all
sufficiently large $R$.

\vskip 5pt \noindent {\bf Proof of  Proposition \ref{prop5}.} Let $n\geq2$,
and $L$ be a differential operator defined by \eqref{(1)}
in $ {\mathbb R}^n$  such that its coefficients  satisfy 
condition \eqref{(14)} for all sufficiently large $R$,  with  some constants $A>0$ and $n-2>\sigma > -2$. Then,  for any  $q>1$ and
any $\nu \in (0,1)\cap (0,q-1)$,  taking $p=2(q-\nu)/(q-1)$ in \eqref{(15)} we
obtain the inequality
\begin{eqnarray}\label{(110)}
{\mathrm {cap}}_{L,p} (\overline {\mathcal B}_{R/2},{\mathbb
R}^n\setminus {\mathcal B}_R ;{\mathbb R}^n)\leq {\hat
C}R^{(n-\sigma-2)(q- (n-\nu(\sigma+2))/({n-\sigma-2}))/(q-1)},
\end{eqnarray}
which holds for all sufficiently large $R$, with some
positive constant $\hat C$  which depends only on $A$, $n$, $p$  and $\sigma$. From  \eqref{(110)}  we obtain that for any $\sigma < n-2$ and any
$1<q\leq (n-\nu(\sigma+2))/({n-\sigma-2})$ the $(L, p)$-capacity of
the condenser $({\overline {\mathcal B}}_{R/2},{\mathbb
R}^n\setminus {\mathcal B}_R ;{\mathbb R}^n)$ is bounded above by a
constant which depends only on $A$, $n$, $p$ and $\sigma$, for all
sufficiently large $R$.

\vskip 5pt \noindent {\bf Proof of Proposition \ref{prop6}.} Let $n\geq2$, and
$L$ be a differential operator defined by \eqref{(1)} in $
{\mathbb R}^n$  such that its coefficients  satisfy 
condition \eqref{(14)} for all sufficiently large $R$,  with  some constants $A>0$ and
$n-2>\sigma>-2$. Then, from \eqref{(15)}  we obtain  the
inequality
\begin{eqnarray}\label{(111)}
{\frak C}_{L, p_1,p_2}(R) \leq c R^{(2n- p_1(\sigma+2))/4}R^{(2n-
p_2(\sigma+2))/2p_2}
\end{eqnarray}
for all sufficiently large $R$, with some positive constant $c$ 
which depends only on $A$, $n$, $p_1$, $p_2$  and $\sigma$. Hence
 for any $\nu \in (0,1)\cap (0,(\sigma+2)/ (n-\sigma -2))$, 
choosing $p_1\geq 2(n-\nu(n-\sigma-2))/(\sigma+2)$ and $p_2\geq
2n/(\sigma+2-\nu(n-\sigma-2))$ in \eqref{(111)} we obtain that the quantity $ {\frak
C}_{L, p_1,p_2}(R)$ is bounded by a positive constant, which depends
only on $A$, $n$, $p_1$, $p_2$ and $\sigma$, for all sufficiently
large $R$, and this concludes the proof of Proposition \ref{prop6}.

\vskip 5pt \noindent {\bf Proof of Proposition \ref{prop7}.} Let $n\geq 2$,
and $L$ be a differential operator defined by \eqref{(1)}
in $ {\mathbb R}^n$ such that its coefficients  satisfy 
condition \eqref{(14)} for all sufficiently large $R$,  with  some constants $A>0$ and
$n-2>\sigma> -2$. Then,  for any $q>1$ and
any $\nu \in (0,1)\cap (0,q-1)$,  letting $p_1=2(q-\nu)/(q-1)$ and
$p_2=2q/(q-1-\nu)$ in \eqref{(111)},
 we obtain the inequality
\begin{eqnarray}\label{(112)}
{\frak C}_{L, p_1,p_2}(R) \leq c
R^{(2q-1-\nu)(n-\sigma-2)(q-n/(n-\sigma-2))/(2q(q-1))}
\end{eqnarray}
for all sufficiently large $R$, with some positive constant $c$ 
which depends only on $A$, $n$, $p_1$, $p_2$  and $\sigma$. From \eqref{(112)},  we have that for any $q\leq n/(n-\sigma-2)$ the
quantity ${\frak C}_{L, p_1,p_2}(R)$ is bounded above  by a
constant,  which depends only on $A$, $n$, $p_1$, $p_2$  and
$\sigma$, for all sufficiently large $R$, and this concludes the
proof of Proposition \ref{prop7}.

\vskip 5pt \noindent {\bf Proof of Proposition \ref{prop8}.} Let $n\geq2$, and
 $L$ be a differential operator defined by the relation \eqref{(1)} in $
{\mathbb R}^n$  such that its coefficients satisfy condition
\eqref{(14)} for  all sufficiently large $R$,  with some  constants $A>0$ and $\sigma > -2$. Then letting $p=2$ in \eqref{(15)} and multiplying  both
sides of this inequality by $R^{-n}$ we obtain the inequality
\begin{eqnarray}\label{(113)}
{\mathrm {cap}}_{L,2} (\overline {\mathcal B}_{R/2},{\mathbb
R}^n\setminus {\mathcal B}_R ;{\mathbb R}^n)R^{-n}\leq \hat
CR^{-\sigma-2},
\end{eqnarray}
which holds for all sufficiently large $R$, with some
positive constant $\hat C$  which depends only on $A$, $n$ and $\sigma$. From
\eqref{(113)},  for any  $\sigma>-2$, we obtain equality \eqref{(33}, and this
concludes the proof of Proposition \ref{prop8}.

\newpage

\noindent \textbf{Author's address:}

\vspace{10 pt}

\noindent Vasilii V. Kurta

\noindent Mathematical Reviews

\noindent 416 Fourth Street, P.O. Box 8604

\noindent Ann Arbor, Michigan 48107-8604, USA

\noindent \textbf {e-mail:} vkurta@umich.edu, vvk@ams.org
\end{document}